\documentclass[a4paper]{article}
\usepackage{amsmath,amssymb,amsfonts}
\usepackage{colordvi}
\usepackage[all]{xy}

\def\mpoint{\;.}
\def\mvirg{\;,}

\def\mpn{\medskip\par\noindent}

\def\mmpn{\vskip 1em minus 1em\par\noindent}

\def\smp{\smallskip\par}

\def\Id{\hbox{\rm id}}
\def\Res{\hbox{\rm Res}}
\def\Ind{\hbox{\rm Ind}}
\def\Hom{\hbox{\rm Hom}}

\def\Inf{\hbox{\rm Inf}}
\def\Def{\hbox{\rm Def}}
\def\Ten{\hbox{\rm Ten}}
\def\Iso{\hbox{\rm Iso}}
\def\Conj{\hbox{\rm Conj}}

\def\Ker{\hbox{\rm Ker}}

\def\Defres{\hbox{\rm Defres}}
\def\Teninf{\hbox{\rm Teninf}}
\def\Indinf{\hbox{\rm Indinf}}

\def\normal{\mathop{\underline\triangleleft}}
\def\op{^{op}}
\def\bs{\backslash}
\def\ls#1#2{{\,^{#1}\!#2}}

\def\Z{\mathbb{Z}}

\def\Q{\mathbb{Q}}

\def\F{\mathbb{F}}

\newcommand{\flh}[2]{\mathop{\hbox to  4ex{\rightarrowfill}}_{#2}^{#1}
\limits}

\newcommand{\sumb}[2]{\sum_{{\scriptstyle #1}\atop {\scriptstyle #2}}}

\def\pf{\par\bigskip\noindent{\bf Proof. }}
\def\endpf{\nolinebreak~\leaders\hbox to 1em{\hss\
\hss}\hfill~\raisebox{.5ex}
{\framebox[1ex]{}}\par\bigskip}

\newcommand{\carre}[8]{\begin{array}{ccc}
#1&\mathop{\hbox to  12mm{\rightarrowfill}}^{\displaystyle{#2}}
\limits&#3\\
\llap{$\displaystyle{#4}$}\left\downarrow\vbox to 6mm{}\right. & &
\left\downarrow\vbox to 6mm{}\right.\rlap{$\displaystyle{#5}$}\\
\ & &\\#6&\mathop{\hbox to 12mm{\rightarrowfill}}_{\displaystyle  #7}
\limits&#8\\
\end{array}}

\newcommand{\limind}[1]{\mathop{\lim}_{\displaystyle\longrightarrow
\atop  \scriptstyle{#1}}\limits}
\newcommand{\limproj}[1]{\mathop{\lim}_{\displaystyle\longleftarrow
\atop  \scriptstyle{#1}}\limits}

\makeatletter
\renewenvironment{enumerate}{\ifnum \@enumdepth >3 \@toodeep\else
       \advance\@enumdepth \@ne
       \edef\@enumctr{enum\romannumeral\the\@enumdepth}\list
       {\csname  label\@enumctr\endcsname}{\setlength{\topsep}{1ex}
\setlength{\itemsep}{0 pt}\usecounter
         {\@enumctr}\def\makelabel##1{\hss\llap{##1}}}\fi}{\endlist}

\renewenvironment{itemize}{\ifnum \@itemdepth >3 \@toodeep\else
\advance\@itemdepth \@ne
\edef\@itemitem{labelitem\romannumeral\the\@itemdepth}
\list{\csname\@itemitem\endcsname}{\setlength{\topsep}{1ex}\setlength
{\itemsep}{0pt}\def\makelabel##1{\hss\llap{##1}}}\fi}
{\endlist}

\def\@sect#1#2#3#4#5#6[#7]#8{\ifnum #2>\c@secnumdepth
      \let\@svsec\@empty\else
      \refstepcounter{#1}\edef\@svsec{\csname the#1\endcsname .
\hskip  .5em}\fi
      \@tempskipa #5\relax
       \ifdim \@tempskipa>\z@
         \begingroup #6\relax
           \@hangfrom{\hskip #3\relax\@svsec}{\interlinepenalty \@M
#8\par}
         \endgroup
        \csname #1mark\endcsname{#7}\addcontentsline
          {toc}{#1}{\ifnum #2>\c@secnumdepth \else
                       \protect\numberline{\csname the#1\endcsname}\fi
                     #7}\else
         \def\@svsechd{#6\hskip #3\relax
                    \@svsec #8\csname #1mark\endcsname
                       {#7}\addcontentsline
                            {toc}{#1}{\ifnum #2>\c@secnumdepth \else
                              \protect\numberline{\csname  the#1
\endcsname}\fi
                        #7}}\fi
      \@xsect{#5}}

\def\section{\pagebreak[3]\setcounter{prop}{0}\setcounter{equation}{0}
\@startsection{section}{1}{\z@}{6ex plus  9ex}{3ex}{\center\reset@font
\large\bf}}
\def\subsection{\pagebreak[3]\refstepcounter{prop}\@startsection
{subsection}{2}{\z@}{4ex plus 6ex}{-1em}{\reset@font\bf}}
\def\subsubsection{\@startsection{subsubsection}{3}{\z@}{4ex plus
6ex}{-1em}{\reset@font\it}}
\makeatother

\def\theprop{\thesection.\arabic{prop}}

\renewenvironment{equation}{\refstepcounter{subsection}\refstepcounter
{prop}$$}{\leqno{\bf (\theprop)}$$}

\newenvironment{enonce}[1]{\pagebreak[3]\refstepcounter{prop}\mmpn
{{\bf  \thesection.\arabic{prop}.\ #1.}}\begin{it} }{\end{it}\smp}
\def\thesection{\arabic{section}}

\newcommand{\result}[1]{\begin{enonce}{#1}}
\newcommand{\fresult}{\end{enonce}}

\begin{document}
\centerline{\large\bf A sectional characterization of the Dade group}
\vspace{.5cm}
\centerline{\bf Serge Bouc and Jacques Th\'evenaz}
\vspace{1cm}\par
\begin{footnotesize}
{\bf Abstract :}
Let $k$ be a field of characteristic $p$ , let $P$ be a finite $p$-
group, where $p$ is an odd prime, and let $D(P)$ be the Dade group of
endo-permutation $kP$-modules. It is known that $D(P)$ is detected
via deflation--restriction by the family of all sections of~$P$ which
are elementary abelian of rank~$\leq2$.
In this paper, we improve this result by characterizing $D(P)$ as the
limit (with respect to deflation--restriction maps and conjugation
maps) of all groups $D(T/S)$ where $T/S$ runs through all sections
of~$P$ which are either elementary abelian of rank~$\leq3$ or
extraspecial of order~$p^3$ and exponent~$p$.

\vspace{1ex}\par
\vspace{-.5ex}\par
\end{footnotesize}

\section{Introduction}

\noindent
Endo-permutation modules for finite $p$-groups play an important role
in the representation theory of finite groups and were classified
recently in~\cite{boclassif}. The set of equivalence classes of such
modules is an abelian group $D(P)$ (with respect to tensor product),
called the Dade group of~$P$. An important ingredient for the
classification is a detection theorem, proved in~\cite{cath2}, which
asserts (for odd~$p$) that the product of all deflation--restriction
maps
$$\prod_{(T,S)} \Defres_{T/S}^P : \; D(P) \longrightarrow \prod_
{(T,S)} D(T/S)$$
is injective, where $T/S$ runs through all sections of~$P$
which are elementary abelian of rank~$\leq2$.
The purpose of this paper is to improve this result and characterize
the image of the injective map above, but actually by changing
slightly the target.\par
In order to motivate our result, let us first mention two classical
cases where a similar situation occurs. The first instance is group
cohomology $H^*(G,k)$ where $G$ is a finite group anf $k$ is a field
of characteristic~$p$. An easy detection result asserts that the
restriction map
$H^*(G,k)\to H^*(P,k)$ is injective where $P$ is a Sylow $p$-subgroup
of~$G$.
This is then improved by characterizing the image of $H^*(G,k)$ as
the set of $G$-stable elements in~$H^*(P,k)$.
This improvement can be stated in a more sophisticated way~:
$H^*(G,k)$ is isomorphic to the limit $\limproj{Q} H^*(Q,k)$ where $Q$
runs
through all $p$-subgroups of~$G$ and the limit is taken with respect
to restrictions and conjugations.\par
The second classical case is the ordinary representation ring $R(G)$
of a finite group~$G$. This is detected on restriction to cyclic
subgroups, but in order to characterize the image, Brauer had to
introduce the larger class of Brauer--elementary subgroups. Brauer's
theorem can be stated as follows~: $R(G)$ is isomorphic to the limit
$\limproj{Q} R(Q)$ where $Q$ runs through all Brauer--elementary
subgroups of~$G$ and the limit is taken with respect to restrictions
and conjugations.\par
We are facing a similar situation with the Dade group, except that
the detection map is not a product of restriction maps to subgroups,
but a product of deflation--restriction maps to sections of the group.
When $p$ is odd, the detection family consists of elementary abelian
$p$-groups of rank~$\leq2$. However, as in the case of the ordinary
representation ring, we need to enlarge this family to impose just
the right conditions for the limit. It turns out that we have to
include also elementary abelian $p$-groups of rank~$3$ and
extraspecial groups of order~$p^3$ and exponent~$p$.
Our main result is the following.

\result{Theorem} \label{main} Let $p$ be an odd prime number and
$P$ a finite $p$-group.
Then the natural map
$$
\prod_{(T,S)} \Defres_{T/S}^P : \; D(P) \longrightarrow \limproj
{(T,S)} D(T/S)
$$
is an isomorphism, where $T/S$ runs through all sections of~$P$ which
are either elementary abelian $p$-groups of rank~$\leq3$ or
isomorphic to the extraspecial group of order~$p^3$ and exponent~$p$.
\fresult

As in earlier work (\cite{both1}, \cite{boclassif}), the Dade group
is viewed as a functor on $p$-groups with morphisms being
compositions of restrictions, tensor inductions, deflations,
inflations, and isomorphisms. In more technical terms, $D(-)$ is
viewed as a biset functor and this plays a role again here.\par

More precisely, whenever $\mathcal Y$ is a class of $p$-groups
closed under taking sections, one can consider biset functors defined
only on~$\mathcal Y$, with values in abelian groups.
Let ${\mathcal F}_{\mathcal Y}$ be the category of all such functors
(an abelian category) and write simply $\mathcal F$ for the category
of functors defined on all finite $p$-groups. There is an obvious
forgetful functor
$${\mathcal O}_{\mathcal Y} : {\mathcal F} \to {\mathcal F}_{\mathcal
Y}$$
and we consider its left and right adjoints
$${\mathcal L}_{\mathcal Y} : {\mathcal F}_{\mathcal Y} \to {\mathcal
F} \quad\text{and}\quad
{\mathcal R}_{\mathcal Y} : {\mathcal F}_{\mathcal Y} \to {\mathcal
F} \;.$$
The connection with limits is the following.

\result{Theorem}
With the notation above, for any functor $M\in{\mathcal F}_{\mathcal
Y}$,  we have :
$$
{\mathcal L}_{\mathcal Y}M(P) \cong \limind{(T,S)\in{\mathcal Y}(P)} M
(T/S)
\qquad\text{and}\qquad
{\mathcal R}_{\mathcal Y}M(P) \cong \limproj{(T,S)\in{\mathcal Y}(P)}
M(T/S) \,,
$$
where ${\mathcal Y}(P)$ denotes the set of pairs $(T,S)$ of subgroups
of~$P$ such that
$S\normal T$ and $T/S\in\nolinebreak{\mathcal Y}$.
\fresult
This type of result is more or less standard in category theory. In
our case, it is made explicit in an appendix at the end of the paper.
The main point is that the various inverse limits can be organized to
yield a functor, namely ${\mathcal R}_{\mathcal Y}M$. Therefore we
can use the machinery of functors throughout this paper.\par
Using this point of view, the natural map of Theorem~\ref{main}
becomes the unit morphism of the adjunction and our main result takes
the following form.
\result{Theorem} Let $p$ be an odd prime number and let $D$ be the
Dade functor. Let ${\mathcal X}_3$ be the class of all $p$-groups
which are either elementary abelian $p$-groups of rank~$\leq3$
or isomorphic to the extraspecial group of order~$p^3$ and exponent~$p$.
Then the unit morphism
$$
\eta_D^{\mathcal{X}_3} : D\to \mathcal{R}_{\mathcal{X}_3}\mathcal{O}_
{\mathcal{X}_3}D
$$
is an isomorphism.
\fresult

Let us now describe the main ingredients of the proof. We consider
first the class ${\mathcal X}$ of all $p$-groups which are either
elementary abelian $p$-groups (without condition on the rank) or
extraspecial groups of order~$p^3$ and exponent~$p$. We prove the
theorem for $\mathcal X$ instead of~${\mathcal X}_3$ and then we use
a technical lemma about elementary abelian groups to derive the more
precise result for~${\mathcal X}_3$.\par
Let $D_{tors}$ be the torsion subfunctor of~$D$.
It follows from \cite{boburn} and \cite{boclassif} that
$D/D_{tors}$ is isomorphic to the $\Z$-dual~$K^*$ of the
functor~$K=\Ker(B\to R_{\Q})$, where $B$ is the Burnside functor
and $R_{\Q}$ is the rational representation functor.
Thus we have an exact sequence
$$0\; \flh{}{} \, D_{tors} \; \flh{}{} \, D \; \flh{}{} \, K^* \; \flh
{}{} \, 0$$
and we consider seperately $D_{tors}$ and $K^*$.\par

For $D_{tors}$, a detection theorem from~\cite{cath2} together with a
result from~\cite{both2} about elementary abelian groups imply that
the unit morphism
$$
\eta_{D_{tors}}^{\mathcal{X}} : D_{tors} \longrightarrow
\mathcal{R}_{\mathcal{X}}\mathcal{O}_{\mathcal{X}}D_{tors}
$$
is an isomorphism. Turning now to~$K^*$, we use an induction theorem
for~$K$ which is proved in~\cite{boclassif} (and which plays a
crucial role for the classification of endo-permutation modules).
This theorem implies that the counit map
$\mathcal{L}_{\mathcal{X}}\mathcal{O}_{\mathcal{X}}K \to K$ is
surjective (and it is here that the extraspecial group of order~$p^3$
is necessary). Dualizing and using again the result from~\cite
{both2}, we obtain that the unit morphism
$$
\eta_{K^*}^{\mathcal{X}} : K^* \longrightarrow
\mathcal{R}_{\mathcal{X}}\mathcal{O}_{\mathcal{X}} K^*
$$
is an isomorphism. The main theorem follows from the two isomorphisms.
\par

Our methods do not work when $p=2$ for several reasons. First $D$ is
not a biset functor if $p=2$. Moreover, the detection map to
elementary abelian sections is not injective in general when $p=2$
and one needs to add the cyclic group of order 4 and
the quaternion group of order~8 in the detecting family.
These two problems are not essential and the arguments could
probably be modified accordingly, but there is one more problem~:
The result from~\cite{both2} concerning elementary abelian
groups is about $p$-torsion and is used for the 2-torsion
group~$D_{tors}$. This collapses when $p=2$. \par

Finally, let us mention that the main result of this article plays a
crucial role in a forthcoming paper of the first author~\cite{boglue}.

\section{Biset functors}

\noindent
Biset functors play a key role throuhout this paper.
We collect in this section the material needed later concerning those
functors.
In particular, we revisit an induction theorem for one specific functor,
called~$K$, which already was important for the classification
of endo-permutation modules in~\cite{boclassif} and
which is used again in an essential way in the present paper.\par

If $Q$ and $P$ are $p$-groups, a $(Q,P)$-biset is a finite set $U$
with a left action of~$Q$ and a right action of~$P$, such that
$(x \cdot u)\cdot y = x\cdot (u\cdot y)$ for all $x\in Q$, $y\in P$,
$u\in U$.
We let ${\mathcal C}_p$ be the category of all finite $p$-groups with
morphisms defined by~:
$$\Hom_{{\mathcal C}_p}(P,Q)=B(Q\times P^{op})\,,$$
the Burnside group of all finite $(Q,P)$-bisets. The composition of
morphisms is $\Z$-linear and is induced by the usual product of
bisets, namely $V\circ U = V\times_QU$ whenever
$V$ is an $(R,Q)$-biset and $U$ is a $(Q,P)$-biset.
A {\em biset functor\/} is an additive functor from~${\mathcal C}_p$
to the category $\mathcal{A}b$ of abelian
groups (see \cite{bodouble} and~\cite{both1}).\par

We now recall how, for suitable choices of bisets~$U$, we obtain
morphisms of restriction,
induction, inflation, deflation, and isomorphisms.
The word `section' was used informally in the introduction,
but we now give a more precise definition.
A {\em section} of a group $P$ is a pair $(T,S)$ of subgroups
of~$P$ such that $S\normal T$.
The group~$T/S$ will be referred to as a {\em subquotient} of~$P$.
If $F$ is a biset functor and $(T,S)$ is a section of~$P$,
then the following homomorphisms are defined:
\par\noindent $\bullet\;$ $\Inf_{T/S}^T : F(T/S) \to F(T)$, induced
by the $(T,T/S)$-biset $T/S$ (inflation);
\par\noindent $\bullet\;$  $\Ind_T^P: F(T) \to F(P)$, induced
by the $(P,T)$-biset $P$ (induction);
\par\noindent $\bullet\;$ $\Def_{T/S}^T : F(T) \to F(T/S)$, induced
by the $(T/S,T)$-biset $S\bs T$ (deflation);
\par\noindent $\bullet\;$  $\Res_T^P: F(P) \to F(T)$, induced
by the $(T,P)$-biset $P$ (restriction);
\par\noindent $\bullet\;$  Whenever $\alpha:P\to Q$ is an
isomorphism, the corresponding isomorphism $\Iso_\alpha:F(P)\to F(Q)$
is induced by the $(Q,P)$-biset $P$ with left action of~$Q$
via~$\alpha^{-1}$. In particular, $\Conj_x: F(T/S)\to F(\ls xT/\!\ls
xS)$
(where $x\in P$) is induced by the
$(\ls xT/\!\ls xS,T/S)$-biset $T/S$ (conjugation by~$x$).\par

It can be shown that any transitive biset can be decomposed as
a product of the above five types of bisets, namely
a composition of a restriction, a deflation, an isomorphism,
an inflation, and an induction
(see Lemma~7.4 in~\cite{both1} or Lemma 3
in~\cite{bodouble}),
so that a biset functor is essentially a functor
endowed with those five types of morphisms.
We also consider the composites~:
\par\noindent $\bullet\;$  $\Indinf_{T/S}^P=\Ind_T^P\Inf_{T/S}^T$,
induced by the $(P,T/S)$-biset
$P\times_T(T/S) \cong P/S$;
\par\noindent $\bullet\;$ $\Defres_{T/S}^P=\Def_{T/S}^T\Res_T^P$,
induced by the $(T/S,P)$-biset
$(S\bs T)\times_T P \cong S\bs P$.\par

If $F$ is a biset functor and $x\in F(Q)$ for some fixed
$p$-group~$Q$, the subfunctor $F_x$ generated by~$x$
is the intersection of all subfunctors $E$ of~$F$
such that $x\in E(Q)$. Equivalently
$$F_x(P)=\Hom_{{\mathcal C}_p}(Q,P)\times_Q x \mvirg $$
and this means in practice that we have to 
consider linear combinations of elements obtained by applying to $x$ 
successively restrictions, deflations, isomorphisms, inflations,
and inductions in order to 
get all possible elements of~$F_x(P)$.\par

Some standard constructions give rise to biset functors.
This is the case in particular for the Burnside group $B(P)$ and the
rational representation group~$R_{\Q}(P)$, which turn out to have strong
connections with the Dade group, as we shall see in the next
section.\par

Recall that $B(P)$ is the Grothendieck group of finite $P$-sets, with
$\Z$-basis consisting of all transitive $P$-sets $P/S$, where $S$ runs
over all subgroups of~$P$ up to conjugation. For every
$(Q,P)$-biset~$U$,
we have a map $B(U): B(P)\to B(Q)$ defined by $X\mapsto U\times_PX$
for every $P$-set~$X$.
This provides a structure of biset functor on~$B$,
called simply the Burnside functor.
Note that if the morphism~$U$ is a restriction, induction, inflation,
or deflation, then the morphism $B(U)$ is indeed the corresponding
morphism between Burnside
groups.\par
On the other hand $R_{\Q}(P)$ is the Grothendieck group of
all ${\Q}P$-modules (that~is, the rational representations of~$P$) and
the irreducible rational representations form a ${\Z}$-basis
of~$R_{\Q}(P)$.
For every $(Q,P)$-biset~$U$, we have a map $R_{\Q}(P)\to R_{\Q}(Q)$
defined by $V\mapsto {\Q}U\otimes_{{\Q}P}V$ for every
${\Q}P$-module~$V$ (where ${\Q}U$ is the permutation
$({\Q}Q,{\Q}P)$-bimodule
with ${\Q}$-basis~$U$). Thus $R_{\Q}$ is a biset functor and again we
recover the usual morphisms of restriction, induction, inflation, or
deflation for representation 
groups.\par
There is a natural morphism $B\to R_\Q$,
which is surjective for~$p$-groups by a result of Ritter and Segal.
We define
$$ K = \Ker(B\to R_\Q)\,.$$
This is a subfunctor of~$B$ which plays a crucial role in the sequel.
We shall need in particular the following induction theorem for the
biset functor~$K$, proved in Section~6 of~\cite{boclassif}. We also
include some additional results which are only implicit
in~\cite{boclassif}.

\result{Theorem} \label{induction}
Let $p$ be an odd prime, let ${\mathcal E}_2$ be the class of all
elementary abelian $p$-groups of rank~$\leq2$, and let ${\mathcal X}_2$
be the class of $p$-groups consisting of ${\mathcal E}_2$ and the
groups isomorphic to the extraspecial
$p$-group~$X$ of order~$p^3$ and exponent~$p$.
\begin{itemize}
\item[(a)] $K$ is generated by some specific element $\delta\in K(X)$.
In particular, for every $p$-group~$P$, the sum of all
$\Indinf_{T/S}^PK(T/S)$, where $(T,S)\in
\mathcal{X}_2(P)$, is equal to~$K(P)$.
\item[(b)] Let $E$ be an elementary abelian group of order~$p^2$ and
let $\epsilon$ be a generator of $K(E)\cong\Z$.
For every $p$-group~$P$, the sum of all $\Indinf_{T/S}^PK(T/S)$,
where $(T,S)\in\mathcal{E}_2(P)$, is equal to $K_\epsilon(P)$,
where $K_\epsilon$ is the subfunctor of~$K$ generated by~$\epsilon$.
\item[(c)] We have $p\cdot K(P) \subseteq K_\epsilon(P)$,
for every $p$-group~$P$.
\end{itemize}

\fresult

\pf
(a) This is exactly Theorem 6.12 and Corollary 6.16
in~\cite{boclassif}.\par

(b) In the Burnside group $B(E)$, define the element
$$\epsilon = E/1 - \sumb{F\leq E}{|F|=p} E/F + p E/E \,.$$
It is easy to see that the corresponding rational representation is
zero, so that $\epsilon \in K(E)$. More precisely, $K(E)\cong\Z$
generated by~$\epsilon$, because $B(E)$ is $\Z$-free of
rank $p+3$ and $R_\Q(E)$ is $\Z$-free of rank $p+2$
(which is the number of cyclic subgroups of~$E$).\par

Define $K_\epsilon$ to be the subfunctor of~$K$ generated by~$\epsilon$.
On evaluation at a $p$-group~$P$, we have
$$K_\epsilon(P)
= \sum_{(T,S)\in\mathcal{E}_2(P)}  \Indinf_{T/S}^P \, K(T/S)$$
because $K(Q)=0$ for any group~$Q$ of order 1 or~$p$ so that any
proper deflation or restriction of~$\epsilon$ is zero.\par

(c) Since $K$ is generated by~$\delta$, the inclusion $p K
\subseteq K_\epsilon$ will follow if we prove that $p\delta$ belongs
to~$K_\epsilon(X)$.
The expression for~$\delta$ in the Burnside group~$B(X)$ is the
following (see~6.9 of~\cite{boclassif})~:
$$\delta = X/I-X/IZ-X/J+X/JZ \,,$$
where $I$ and $J$ are non-conjugate subgroups of order~$p$ and $Z$ is
the centre of~$X$. Now $IZ$ is elementary abelian of order~$p^2$ and
the subgroup~$I$ has $p$ conjugates in~$X$ contained in~$IZ$ (the
only additional subgroup of order~$p$ in~$IZ$ being~$Z$). Therefore,
if $\epsilon_{IZ}$ denotes~$\epsilon$ viewed in~$B(IZ)$, we have
$$\Ind_{IZ}^X \, \epsilon_{IZ} = X/1 - X/Z - pX/I + pX/IZ
$$
and similarly when $I$ is replaced by~$J$. It follows that
$$\Ind_{JZ}^X \, \epsilon_{JZ} - \Ind_{IZ}^X \, \epsilon_{IZ} = p\delta
$$
and this shows that $p\delta \in K_\epsilon(X)$.
\endpf

We end this section by considering duality (defined in~\cite{boburn}).
The structure of biset functor on a dual uses opposite bisets as
follows.
If $F$ is a biset functor, define $F^*(P)=\Hom_\Z(F(P),\Z)$
for every $p$-group~$P$.
If $U$ is a $(Q,P)$-biset, then define $F^*(U):F^*(P)\to F^*(Q)$ to be
the transpose
$$F^*(U) = F(U\op)^t \mvirg$$
where $U\op$ denotes the $(P,Q)$-biset with underlying set equal to~$U$
and left $P$-action of~$x\in P$ defined by using the given right action
of~$x^{-1}$ on~$U$, and similarly for the action of~$Q$.
The effect of this definition is that restriction for~$F^*$ is the
transpose of induction for~$F$ (and conversely), while inflation
for~$F^*$ is the transpose of deflation for~$F$ (and conversely).\par

By definition of~$K$ and by the result of Ritter and Segal, we have an
exact sequence
$$0\; \flh{}{} \, K \, \flh{}{} B \flh{}{}  R_\Q \flh{}{} \; 0 \mvirg$$
and since we get free $\Z$-modules on each evaluation, we obtain by
duality
an exact sequence
\begin{equation}\label{exact-K*}
0\; \flh{}{} \, R_\Q^* \, \flh{}{} B^* \flh{}{}  K^* \flh{}{} \; 0
\mpoint
\end{equation}

This exact sequence has an intimate connection with the Dade functor,
as we shall recall in the next section.

\section{The Dade functor}

\noindent
We collect in this section the various results we need about the Dade
group and the Dade functor.
Let $p$ be a prime number (which will soon be assumed to be odd) and
let $k$ be a field of characteristic~$p$. For any finite $p$-group~$P$,
we let $D(P)$ be the Dade group of~$P$, that is, the group of
equivalence classes of endo-permutation $kP$-modules (see~\cite
{both1} for details). If $X$ is a 
non empty finite $P$-set, the kernel of the
augmentation map $kX\to k$ is called a relative syzygy of the trivial
module (relative to~$X$) and is an endo-permutation module.
These are the main examples of endo-permutation modules and
we let $D^\Omega(P)$ denote the subgroup of~$D(P)$ generated
by relative syzygies of the trivial module.\par

Recall that if $(T,S)$ is a section of~$P$, we have a restriction map
$\Res_T^P: D(P) \to D(T)$ and a tensor induction map
$\Ten_T^P: D(T)  \to D(P)$, as well as an
inflation map $\Inf_{T/S}^{T}: D(T/S) \to D(T)$ and a deflation map
$\Def_{T/S}^{T}: D(T) \to D(T/S)$.
We also write  $\Defres_{T/S}^P = \Def_{T/S}^T \, \Res_T^P$ and
similarly
$\Teninf_{T/S}^P =\Ten_T^P \, \Inf_{T/S}^T$.
It is proved in~\cite{bomega} that $D^\Omega$ is invariant under those
operations.\par

The main detection theorem which we are concerned with is the
following (Theorem~13.1 in~\cite{cath2}).

\result{Theorem} \label{detection}
Let ${\mathcal E}_2(P)$ be the set of sections $(T,S)$ of~$P$ such
that $T/S$ is elementary abelian of rank~$\leq2$.
If $p$ is odd, the map
$$\prod_{(T,S)\in\mathcal{E}_2(P)} \Defres_{T/S}^P~:
D(P) \longrightarrow \prod_{(T,S)\in\mathcal{E}_2(P)} D(T/S) $$
is injective.
\fresult

When $p$ is odd, we want to view $D$ as a biset functor.
Any $(Q,P)$-biset~$U$ induces a group homomorphism
$D(U)~: D(P)\to D(Q)$ (see Corollary~2.13  of~\cite{both1}).
For the specific choices of bisets~$U$ described in Section~2,
we recover the above morphisms of restriction, tensor
induction, inflation, deflation, and isomorphisms
(see Section~2 of~\cite{both1}).\par

\result{Proposition}
\begin{itemize}
\item[(a)] $D^\Omega$ is a biset functor.
\item[(b)] If $p$ is odd, $D$ is a biset functor.
Moreover, $D_{tors}$ is a subfunctor of~$D$, where
$D_{tors}(P)$ denotes the torsion subgroup of~$D(P)$.
\end{itemize}
\fresult

\pf (a) The proof that $D^\Omega$
is a biset functor is the starting point of~\cite{boburn}.

(b) There are two ways of proving that $D$ is a biset functor.
If $p$ is odd, we have $D^\Omega(P)=D(P)$ by Theorem~7.7
of~\cite{boclassif}, which is one of the main theorems of the
classification
of endo-permutation $kP$-modules.
Thus $D=D^\Omega$ and $D$ is a biset functor.\par

The other proof does not require the full classification but uses
Theorem~\ref{detection}.
This approach is explicit in Theorem~10.1 of~\cite{survey}.\par

The fact that $D_{tors}$ is a subfunctor is clear since a torsion
element is mapped to a torsion element by a group homomorphism.
\endpf

Now we recall the connection between $D$ and the standard biset
functors considered in Section~2.

\result{Theorem} \label{exact}
Let $p$ be an odd prime.
\begin{itemize}
\item[(a)] The natural injection $R_\Q^* \to B^*$ fits into an exact
sequence of biset functors
$$0\; \flh{}{} \, R_\Q^* \, \flh{}{} B^* \flh{}{}  D^\Omega/D^\Omega_
{tors} \flh{}{} \; 0 \mpoint$$
\item[(b)]
There is an exact sequence of biset functors
$$0\; \flh{}{} \, D_{tors} \, \flh{}{} D \flh{}{}  K^* \flh{}{} \; 0
\mpoint$$
\end{itemize}
\fresult

\pf (a) This exactly Theorem 1.8 of~\cite{boburn}.\par

(b) Since $p$ is odd, we have $D^\Omega=D$ by Theorem~7.7
of~\cite{boclassif}. Comparing the exact sequence of~(a) with
the sequence~\ref{exact-K*},
we see that $K^* \cong D/D_{tors}$, hence the result.
\endpf

\section{Limits and units of adjunctions}

\noindent
Let $\mathcal Y$ be a class of finite $p$-groups.
We shall say that $\mathcal Y$ is {\em closed under taking sections\/}
if for any $R\in\mathcal Y$ and any section $(T,S)$ of~$R$,
any group isomorphic to~$T/S$ belongs to~$\mathcal Y$.
In particular $\mathcal Y$ is closed under
isomorphisms. For every $p$-group~$P$, we define ${\mathcal Y}(P)$ to
be the set of all sections $(T,S)$ of~$P$ such that $T/S\in \mathcal
Y$.\par

Let $\mathcal Y$ be a class of finite $p$-groups, closed under taking
sections.
Let ${\mathcal C}_{\mathcal Y}$ be the full subcategory of the
category~${\mathcal C}_p$ whose objects are in~$\mathcal Y$.
A biset functor on~$\mathcal Y$ is an additive functor
from ${\mathcal C}_{\mathcal Y}$ to the category of abelian groups.
Let ${\mathcal F}$ be the category of biset functors (defined on the
whole of~${\mathcal C}_p$) and
let ${\mathcal F}_{\mathcal Y}$ be the category of biset functors
on~$\mathcal Y$. Both are abelian categories. There is an obvious
forgetful functor
$${\mathcal O}_{\mathcal Y} : {\mathcal F} \longrightarrow {\mathcal
F}_{\mathcal Y}$$
and we wish to construct functors in the opposite direction. If $M$
is a biset functor on~$\mathcal Y$, we define, for every $p$-group~$P$,
$${\mathcal L}_{\mathcal Y}M(P) \cong \limind{(T,S)\in{\mathcal Y}
(P)} M(T/S)
\qquad\text{and}\qquad
{\mathcal R}_{\mathcal Y}M(P) \cong \limproj{(T,S)\in{\mathcal Y}(P)}
M(T/S) \,.
$$
In the first case, the direct limit is constructed with respect to
all induction--inflation maps and all conjugation maps. In the second
case, the inverse limit is constructed with respect to all
deflation--restriction maps and all conjugation maps
(see~(\ref{lim}) for details).
For practical purposes, we always view the inverse limit
$\limproj{(T,S)\in{\mathcal Y}(P)} M(T/S)$ as a subset of the direct
product
$\prod_{(T,S)\in{\mathcal Y}(P)} M(T/S)$.\par

It turns out that it is possible to put a structure of biset functors
on~${\mathcal L}_{\mathcal Y}M$ and~${\mathcal R}_{\mathcal Y}M$.
More precisely, we have the following theorem.

\result{Theorem} \label{adjoints}
Let $\mathcal Y$ be a class of finite $p$-groups, closed under taking
sections.
\begin{itemize}
\item[(a)] For every biset functor $M$ defined on~$\mathcal Y$, there
is a biset functor ${\mathcal L}_{\mathcal Y}M$ defined on all
$p$-groups such that, for every $p$-group~$P$,
the evaluation of~${\mathcal L}_{\mathcal Y}M$ at~$P$ is given by the
direct limit above.
Moreover, the assignment $M \mapsto {\mathcal L}_{\mathcal Y}M$
defines a functor
${\mathcal F}_{\mathcal Y} \longrightarrow {\mathcal F}$ which is
left adjoint to the forgetful functor~${\mathcal O}_{\mathcal Y}$.
\item[(b)] For every biset functor $M$ defined on~$\mathcal Y$, there
is a biset functor ${\mathcal R}_{\mathcal Y}M$ defined on all
$p$-groups such that, for every $p$-group~$P$,
the evaluation of~${\mathcal R}_{\mathcal Y}M$ at~$P$ is given by the
inverse limit above.
Moreover, the assignment $M \mapsto {\mathcal R}_{\mathcal Y}M$
defines a functor
${\mathcal F}_{\mathcal Y} \longrightarrow {\mathcal F}$ which is
right adjoint to the forgetful functor~${\mathcal O}_{\mathcal Y}$.
\end{itemize}
\fresult

The proof is rather technical and is independent of the main purpose
of this paper. It is included in the appendix.
Because of the adjunction, there is a unit natural transformation
$$\eta^{\mathcal Y} : {\rm id}_{\mathcal F} \longrightarrow
{\mathcal R}_{\mathcal Y}{\mathcal O}_{\mathcal Y} \,,$$
and we write its evaluation at a biset functor~$F$ as follows~:
$$\eta^{\mathcal Y}_F : F \longrightarrow
{\mathcal R}_{\mathcal Y}{\mathcal O}_{\mathcal Y}F \,.$$
It is also proved in the appendix (Corollary~\ref{unit-defres}) that
when $\eta^{\mathcal Y}_F$ is evaluated at a
\linebreak[2]$p$-group~$P$, one gets
a map which coincides with the product of deflation--restrictions
to~${\mathcal Y}(P)$~:
$$\eta^{\mathcal Y}_{F}(P) = \prod_{(T,S)\in{\mathcal Y}(P)} \!\!
\Defres^P_{T/S} :
\; F(P) \longrightarrow
{\mathcal R}_{\mathcal Y}{\mathcal O}_{\mathcal Y}F(P) =  \limproj
{(T,S)\in{\mathcal Y}(P)} F(T/S)\,.$$
The discussion of this map is at the heart of the present paper.\par

We shall need to compare the constructions corresponding to two
different classes of $p$-groups
$\mathcal Y$ and $\mathcal Z$ such that ${\mathcal Z}\subseteq
{\mathcal Y}$.
In that case, for every $p$-group~$P$, there is a natural homomorphism
$$\pi_{\mathcal{Z}}^{\mathcal Y}(P) : {\mathcal R}_{\mathcal Y}
{\mathcal O}_{\mathcal Y}F(P)
\longrightarrow
{\mathcal R}_{\mathcal Z}{\mathcal O}_{\mathcal Z}F(P)$$
which is simply the restriction of the projection
$$\prod_{(T,S)\in{\mathcal Y}(P)} F(T/S) \longrightarrow
\prod_{(T,S)\in{\mathcal Z}(P)} F(T/S) \,.$$
It can be shown that this defines a morphism of functors
$\pi_{\mathcal Z}^{\mathcal Y} : {\mathcal R}_{\mathcal Y}{\mathcal O}
_{\mathcal Y}F
\longrightarrow
{\mathcal R}_{\mathcal Z}{\mathcal O}_{\mathcal Z}F$,
but we actually only need here
the evaluation~$\pi_{\mathcal Z}^{\mathcal Y}(P)$ at~$P$.

Our first lemma is concerned with the question of passing from a
class~$\mathcal Z$ to a larger class~$\mathcal Y$.

\result{Lemma} \label{agrandir}
Let $\mathcal Y$ and $\mathcal Z$ be two classes of $p$-groups,
closed under taking sections, such that ${\mathcal Z}\subseteq
{\mathcal Y}$.  Let $F$ be a biset functor on $p$-groups.
\begin{itemize}
\item[(a)] If the unit morphism $\eta^{\mathcal Z}_F : F \to
{\mathcal R}_{\mathcal Z}{\mathcal O}_{\mathcal Z}F$ is injective,
then so are the unit morphism
$\eta^{\mathcal Y}_F : F \to
{\mathcal R}_{\mathcal Y}{\mathcal O}_{\mathcal Y}F$
and the morphism $\pi_{\mathcal Z}^{\mathcal Y} : {\mathcal R}_
{\mathcal Y}{\mathcal O}_{\mathcal Y}F
\to {\mathcal R}_{\mathcal Z}{\mathcal O}_{\mathcal Z}F$.
\item[(b)] If the unit morphism $\eta^{\mathcal Z}_F : F \to
{\mathcal R}_{\mathcal Z}{\mathcal O}_{\mathcal Z}F$ is an
isomorphism, then so are the unit morphism
$\eta^{\mathcal Y}_F : F \to
{\mathcal R}_{\mathcal Y}{\mathcal O}_{\mathcal Y}F$
and the morphism $\pi_{\mathcal Z}^{\mathcal Y} : {\mathcal R}_
{\mathcal Y}{\mathcal O}_{\mathcal Y}F
\to {\mathcal R}_{\mathcal Z}{\mathcal O}_{\mathcal Z}F$.
\end{itemize}
\fresult

\pf
(a) Since $\mathcal{Z}\subseteq \mathcal{Y}$, for any $p$-group $P$,
we have a commutative diagram
$$\xymatrix{
F(P)\ar[r]^-{\eta_{F}^\mathcal{Y}(P)}
\ar[rd]_-{\eta_{F}^\mathcal{Z}(P)}
&\mathcal{R}_{\mathcal{Y}}\mathcal{O}_{\mathcal{Y}}F(P)
\ar[d]^{\pi_{\mathcal Z}^{\mathcal Y}(P)}\phantom{\mpoint}\\
&\mathcal{R}_\mathcal{Z}\mathcal{O}_\mathcal{Z}F(P)
}
$$
which shows that $\eta_F^\mathcal{Y}$ is injective if $\eta_F^\mathcal
{Z}$ is.
Moreover if $u\in\Ker\, \pi_{\mathcal Z}^{\mathcal Y}(P)$, view $u$
as a sequence $(u_{T,S})_{(T,S)\in\mathcal{Y}(P)}$, with some
compatibility conditions. Fix $(T,S)$ in $\mathcal{Y}(P)$,
and consider $(T'/S,S'/S)\in \mathcal{Z}(T/S)$, or equivalently
$(T',S')\in\mathcal{Z}(P)$ with ${S\leq S'\leq T'\leq T}$.
Then $\Defres_{T'/S'}^{T/S}u_{T,S}=u_{T',S'}=0$ and this holds for
every section in~${\mathcal Z}(T/S)$.
It follows that $u_{T,S}$ is in the kernel of the map
$$\eta_{F}^\mathcal{Z}(T/S) : F(T/S) \to
{\mathcal R}_{\mathcal Z}{\mathcal O}_{\mathcal Z}F(T/S) \,.$$
Thus $u_{T,S}=0$, since $\eta^{\mathcal Z}_F$ is injective by
assumption.
So $\pi_{\mathcal Z}^{\mathcal Y}(P)$ is injective.
\par

(b) If $\eta_{F}^\mathcal{Z}(P)$ is an isomorphism, then the diagram
above shows that $\pi_{\mathcal Z}^{\mathcal Y}(P)$ is surjective.
Hence by~(a), it is an isomorphism. It now follows that
$\eta_{F}^\mathcal{Y}(P)$ is also an isomorphism.
\endpf

The next lemma is concerned with the question of passing from a class~
$\mathcal Y$
to a smaller class~$\mathcal Z$, provided a suitable assumption holds.

\result{Lemma} \label{rapetisser}
Let $\mathcal Y$ and $\mathcal Z$ be two classes of $p$-groups,
closed under taking sections, such that ${\mathcal Z}\subseteq
{\mathcal Y}$.  Let $F$ be a biset functor on $p$-groups.
If the unit morphism $\eta^{\mathcal Y}_F : F \to
{\mathcal R}_{\mathcal Y}{\mathcal O}_{\mathcal Y}F$ is an
isomorphism and if,
for every $Q\in{\mathcal Y}$, the evaluation at~$Q$ of the unit morphism
$\eta^{\mathcal Z}_F(Q) : F(Q) \to
{\mathcal R}_{\mathcal Z}{\mathcal O}_{\mathcal Z}F(Q)$ is an
isomorphism, then the unit morphism
$\eta^{\mathcal Z}_F : F \to
{\mathcal R}_{\mathcal Z}{\mathcal O}_{\mathcal Z}F$ is an isomorphism.
\fresult

\pf For any $Q\in{\mathcal Y}$, we have a commutative diagram
$$\xymatrix{
F(Q)\ar[r]^-{\eta_{F}^\mathcal{Y}(Q)}
\ar[rd]_-{\eta_{F}^\mathcal{Z}(Q)}
&\mathcal{R}_{\mathcal{Y}}\mathcal{O}_{\mathcal{Y}}F(Q)
\ar[d]^{\pi_{\mathcal Z}^{\mathcal Y}(Q)}\phantom{\mpoint}
&\!\!\!\!\!\!\!\!\!\!\!\!\!\!\!\!\!\!\! = F(Q)\\
&\mathcal{R}_\mathcal{Z}\mathcal{O}_\mathcal{Z}F(Q) &
}
$$
and the top map is the identity.
It follows that $\pi_{\mathcal Z}^{\mathcal Y}(Q)$ can be identified
with $\eta_{F}^\mathcal{Z}(Q)$ and is in particular an isomorphism.
Now for any $p$-group $P$ and any section~$(V,U)$ of~$P$
such that $V/U\in{\mathcal Y}$,
we have a commutative diagram
$$\xymatrix{
F(P)\ar[r]^-{\eta_{F}^\mathcal{Y}(P)}_-{\cong}
\ar[rd]_-{\eta_{F}^\mathcal{Z}(P)}
&\mathcal{R}_{\mathcal{Y}}\mathcal{O}_{\mathcal{Y}}F(P)
\ar[d]^{\pi_{\mathcal Z}^{\mathcal Y}(P)}\phantom{\mpoint}
\ar[rr]^-{\Defres^P_{V/U}}
&&\;\;\mathcal{R}_{\mathcal{Y}}\mathcal{O}_{\mathcal{Y}}F(V/U)
\ar[d]^{\pi_{\mathcal Z}^{\mathcal Y}(V/U)}_\cong\phantom{\mpoint}
& \!\!\!\!\!\!\!\!\!\!\!\!\!\!\!\!\!\!\! = F(V/U) \\
&\mathcal{R}_\mathcal{Z}\mathcal{O}_\mathcal{Z}F(P) \;\;
\ar[rr]^-{\Defres^P_{V/U}}
&& \;\;\mathcal{R}_\mathcal{Z}\mathcal{O}_\mathcal{Z}F(V/U) &
}
$$
and we note that the map $\Defres^P_{V/U}$ on the first line is
actually the projection~on~$F(V/U)$ when the limit $\mathcal{R}_
{\mathcal{Y}}\mathcal{O}_{\mathcal{Y}}F(P)$ is viewed
as a subset of the product $\prod_{(T,S)\in{\mathcal Y}(P)} F(T/S)$
(see Example~\ref{ex1}).\par
When $(V,U)$ runs through the set ${\mathcal Y}(P)$ of all sections
of~$P$ belonging to~$\mathcal{Y}$, the right hand side of the diagram
gives rise to the following commutative square~:
$$\xymatrix{
\mathcal{R}_{\mathcal{Y}}\mathcal{O}_{\mathcal{Y}}F(P)
\ar[rr]^-{\Defres}_-{\cong}
\ar[d]^{\pi_{\mathcal Z}^{\mathcal Y}(P)}\phantom{\mpoint}
&& \raisebox{-4ex}{$\limproj{(V,U)\in{\mathcal Y}(P)}
F(V/U)$}\ar[d]_\cong\phantom{\mpoint} \\
\mathcal{R}_\mathcal{Z}\mathcal{O}_\mathcal{Z}F(P)\;\;
\ar[rr]^-{\Defres}
&& \raisebox{-4ex}{$\limproj{(V,U)\in{\mathcal Y}(P)}
\mathcal{R}_\mathcal{Z}\mathcal{O}_\mathcal{Z}F(V/U)$}
}
$$
We want to prove that $\pi_{\mathcal Z}^{\mathcal Y}(P)$ is an
isomorphism.
The map on the right hand side is an isomorphism by the argument
above. The top map $\Defres$ is an isomorphism by the very definition
of $\mathcal{R}_{\mathcal{Y}}\mathcal{O}_{\mathcal{Y}}F(P)$. Thus it
suffices to prove that the bottom map $\Defres$ is injective.\par

Let $a\in\mathcal{R}_\mathcal{Z}\mathcal{O}_\mathcal{Z}F(P)$ be in
the kernel of~$\Defres$ and write $a=(a_{T,S})$ where $(T,S)$ runs
through $\mathcal{Z}(P)$ and $a_{T,S}\in F(T/S)$.
Since ${\mathcal Z}\subseteq{\mathcal Y}$, we can choose the
component $(V,U)=(T,S)$ in the right hand side limit and we have $
\mathcal{R}_\mathcal{Z}\mathcal{O}_\mathcal{Z}F(T/S)=F(T/S)$, because
$T/S\in\mathcal Z$.
Then the $(T,S)$-component of~$\Defres(a)$ is just~$a_{T,S}$
(see Example~\ref{ex1} applied to the class~$\mathcal Z$).
Since $\Defres(a)=0$, we have $a_{T,S}=0$, that is, $a=0$.
This proves the desired injectivity, hence $\pi_{\mathcal Z}^
{\mathcal Y}(P)$ is an isomorphism.\par

Returning to the previous diagram, we now deduce that $\eta_{F}^
\mathcal{Z}(P)$ is an isomorphism, as was to be shown.
\endpf

For the class of elementary abelian $p$-groups, we have the following
key result.

\result{Theorem} \label{abelem} Let $\mathcal E$ be the class of
elementary abelian $p$-groups.
Let $F$ be a biset functor on $p$-groups and
let $\eta^{\mathcal E}_F : F \to
{\mathcal R}_{\mathcal E}{\mathcal O}_{\mathcal E}F$ be the unit
morphism.
For any fixed $p$-group $P$, there is a group homomorphism
$$\sigma_P : \mathcal{R}_{\mathcal{E}}\mathcal{O}_{\mathcal{E}}F(P)
\to F(P)$$
such that
$\eta_{F}^\mathcal{E}(P) \circ\sigma_P=|P| \cdot  \Id$.
\fresult

\pf This is exactly Theorem~5.1 in~\cite{both2}.
Actually the homomorphism $\sigma_P$ is explicitly defined by the
formula
$$\sigma_P(u)=\sum_{(T,S)\in\mathcal{E}(P)}|S|\mu(S,T)\Indinf_{T/
S}^Pu_{S,\Phi(T)}\mvirg$$
for $u\in \mathcal{R}_{\mathcal{E}}\mathcal{O}_{\mathcal{E}}F(P)$,
where $\mu(S,T)$ is the M\"obius functions of the poset of subgroups
of~$P$.
The proof consists in a computation which shows
that $\eta_{F}^\mathcal{E}(P) \circ\sigma_P=|P| \cdot  \Id$.
\endpf

\result{Corollary} \label{ptorsion} Let $\mathcal Y$ be a class
of $p$-groups, closed under taking sections and containing
the class ${\mathcal E}$ of all elementary abelian $p$-groups.
Let $F$ be a biset functor on $p$-groups such that the unit
morphism $\eta_F^\mathcal{E}$ is injective.
For any fixed $p$-group $P$, there is a group homomorphism
$$\tau_P : \mathcal{R}_{\mathcal{Y}}\mathcal{O}_{\mathcal{Y}}F(P)\to F
(P)$$
such that $\eta_{F}^\mathcal{Y}(P) \circ\tau_P=|P| \cdot  \Id$.
In particular the cokernel of $\eta_{F}^\mathcal{Y}(P)$ is a
$|P|$-torsion group.
\fresult

\pf
Since $\mathcal{Y}\supseteq\mathcal{E}$, we have a commutative diagram
$$\xymatrix{
F(P)\ar[r]^-{\eta_{F}^\mathcal{Y}(P)}
\ar[rd]_-{\eta_{F}^\mathcal{E}(P)}&\mathcal{R}_{\mathcal{Y}}\mathcal
{O}_{\mathcal{Y}}F(P)
\ar[d]^{\pi_{\mathcal E}^{\mathcal Y}(P)}\phantom{\mpoint}\\
&\mathcal{R}_{\mathcal{E}}\mathcal{O}_{\mathcal{E}}F(P)
}
$$
as in the proof of Lemma~\ref{agrandir}.
By Theorem~\ref{abelem}, there exists a map
$$\sigma_P : \mathcal{R}_{\mathcal{E}}\mathcal{O}_{\mathcal{E}}F(P)
\to F(P)$$
such that $\eta_{F}^\mathcal{E}(P)\circ\sigma_P=|P| \cdot  \Id$.
Setting $\tau_P=\sigma_P\circ\pi_{\mathcal E}^{\mathcal Y}(P)$, we have
\begin{eqnarray*}
\pi_{\mathcal E}^{\mathcal Y}(P)\circ \eta_{F}^\mathcal{Y}(P) \circ
\tau_P
&=&\pi_{\mathcal E}^{\mathcal Y}(P)\circ \eta_{F}^\mathcal{Y}(P)
\circ \sigma_P\circ\pi_{\mathcal E}^{\mathcal Y}(P)
=\eta_{F}^\mathcal{E}(P)\circ \sigma_P\circ\pi_{\mathcal E}^{\mathcal
Y}(P) \\
&=&(|P| \cdot \Id) \circ \pi_{\mathcal E}^{\mathcal Y}(P)
= \pi_{\mathcal E}^{\mathcal Y}(P) \circ (|P|\cdot  \Id) \,.
\end{eqnarray*}
But by Lemma~\ref{agrandir}, the map
$\pi_{\mathcal E}^{\mathcal Y}(P)$ is injective,
because we have assumed that the unit morphism
$\eta_F^\mathcal{E}$ is injective.
It follows that $\eta_{F}^\mathcal
{Y}(P)\circ\tau_P=|P|\cdot  \Id$, as was to be shown.
\endpf

\section{The main result}

\noindent
We now return to the Dade functor and prove the main result. We first
work with the class~$\mathcal X$ consisting of all elementary abelian
$p$-groups and extraspecial $p$-groups of order~$p^3$
and exponent~$p$ (where $p$ is odd).
By Theorem~\ref{exact}, there is an exact
sequence of biset functors
$$0\; \flh{}{} \, D_{tors} \, \flh{}{} D \flh{}{}  K^* \flh{}{} \; 0 $$
and we consider $D_{tors}$ and $K^*$ seperately.

\result{Proposition} \label{K*iso} The unit homomorphism
$$\eta_{K^*}^{\mathcal{X}} : K^*\longrightarrow \mathcal{R}_\mathcal
{X}\mathcal{O}_\mathcal{X}K^*$$
is an isomorphism.
\fresult

\pf Consider the counit morphisms
$\varepsilon_{K}^{\mathcal{X}} : \mathcal{L}_\mathcal{X}\mathcal{O}_
\mathcal{X}K\to K$.
By construction of $\mathcal{L}_\mathcal{X}$ and by a result which is
dual to Corollary~\ref{unit-defres},
for any $p$-group $P$, the image of $\varepsilon_{K}^\mathcal{X}(P)$
is equal to the sum of all $\Indinf_{T/S}^PK(T/S)$,
where $(T,S)\in\mathcal{X}(P)$.
Since $\mathcal{X}$ contains~$\mathcal{X}_2$, this sum is the whole
of $K(P)$, by Theorem~\ref{induction}. It is for this crucial result
that we need to include the extraspecial group of order $p^3$ and
exponent~$p$ in our class~$\mathcal X$.
It follows that the morphism $\varepsilon_{K}^{\mathcal{X}}$ is
surjective.\par

Let $M$ denote the kernel of $\varepsilon_{K}^{\mathcal{X}}$, so that
we have a short exact sequence of functors
$$
0\; \flh{}{} \, M \, \flh{}{} \, \mathcal{L}_\mathcal{X}\mathcal{O}_
\mathcal{X}K
\, \flh{}{} \, K \, \flh{}{} \, 0\mpoint
$$
Since $B(P)$ is a free abelian group for any $p$-group $P$, so is
its subgroup~$K(P)$. Therefore, taking $\Z$-duals of this exact
sequence gives the following short exact sequence of biset functors
$$
0 \, \flh{}{} \, K^* \, \flh{}{} \, (\mathcal{L}_\mathcal{X}\mathcal
{O}_\mathcal{X}K)^*
\, \flh{}{} \,  M^* \, \flh{}{} \, 0\mpoint
$$
Now the dual of a colimit is isomorphic to the limit of the duals. So
there is an isomorphism
$$(\mathcal{L}_\mathcal{X}\mathcal{O}_\mathcal{X}K)^*\cong \mathcal{R}
_\mathcal{X}\mathcal{O}_\mathcal{X}K^*\mvirg$$
and the previous sequence becomes
$$
0 \, \flh{}{} \, K^* \, \flh{\eta_{K^*}^{\mathcal{X}}}{} \,  \mathcal
{R}_\mathcal{X}\mathcal{O}_\mathcal{X}K^* \, \flh{}{} \, M^* \, \flh{}
{} \, 0\mpoint
$$
Thus for any $p$-group $P$, we have a short exact sequence of abelian
groups
$$
0 \, \flh{}{} \, K^*(P) \, \flh{\eta_{K^*}^{\mathcal{X}}(P)}{} \,
\mathcal{R}_\mathcal{X}\mathcal{O}_\mathcal{X}K^*(P) \, \flh{}{} \,
M^*(P) \, \flh{}{} \, 0\mpoint
$$
Moreover the $\Z$-dual of any 
abelian group is
torsion free, so $M^*(P)$ is torsion free. We are going to show that
the morphism $\eta_{K^*}^\mathcal{E}$ is injective. It follows that
we can apply Corollary~\ref{ptorsion} and therefore
$M^*(P)$ is a $|P|$-torsion group. Thus $M^*(P)=0$,
and the map $\eta_{K^*}^\mathcal{X}$ is an isomorphism, as was to be
shown.\par

It remains to prove the claim, i.e.~to show that $\eta_{K^*}^\mathcal
{E}(P)$ is injective for any $p$-group~$P$. Let $\varphi\in K^*(P)$
such that $\eta_{K^*}^\mathcal{E}(P)(\varphi)=0$.
By Corollary~\ref{unit-defres}, this means that
$\Defres_{T/S}^P\varphi=0$, for any
section $(T,S)$ of $P$ such that $T/S$ is elementary abelian.
By definition of dual functors (see the end of Section~2), this is
equivalent to
$$\varphi\big
(\Indinf_{T/S}^PK(T/S)\big)=0\mvirg$$
for any $(T,S)\in \mathcal{E}(P)$. By
Theorem~\ref{induction}, we obtain $\varphi(K_\varepsilon(P))=0$.
Hence $p\varphi(K(P))=0$, because $pK(P)\subseteq K_\varepsilon(P)$,
by Theorem~\ref{induction} again.
Since $\varphi$ has values in $\Z$, it follows that $\varphi=0$,
proving the injectivity of~$\eta_{K^*}^\mathcal{E}(P)$.
\endpf

\result{Remark}
{\rm The functor $M$ which appears in the proof is a torsion functor,
because it has the property that $M^*=0$. We do not know if $M=0$ or
if $M$ carries some relevant information.}
\fresult

\bigskip
Now we turn to the analysis of the functor $D_{tors}$.
\result{Proposition} \label{D-tors-iso} Let $\mathcal{Y}$ be a class of $p$-groups, closed under taking subquotients, and containing $\mathcal{E}$. If $p$ is odd, the unit
homomorphism
$$\eta_{D_{tors}}^{\mathcal{Y}} : D_{tors}\longrightarrow \mathcal{R}_
\mathcal{Y}\mathcal{O}_\mathcal{Y}D_{tors}$$
is an isomorphism.

\fresult

\pf By Lemma~\ref{agrandir}, it suffices to consider the case $\mathcal{Y}=\mathcal{E}$. By Theorem~\ref{detection}, the morphism $\eta_{D_{tors}}^\mathcal{E}:
D_{tors}\longrightarrow \mathcal{R}_{\mathcal{E}}\mathcal{O}_{\mathcal{E}}D_{tors}$ is injective.
Moreover, by Corollary~\ref{ptorsion},
for any $p$-group~$P$, the cokernel of
$\eta_{D_{tors}}^\mathcal{E}(P)$ is a $|P|$-torsion group.
But the torsion part of the Dade group of an odd order
$p$-group is a 2-torsion group (this is a consequence of
Theorem~\ref{detection}, see Corollary~13.2 of~\cite{cath2}).
It follows that the cokernel of $\eta_{D_{tors}}^\mathcal{E}(P)$ is
trivial.
Thus $\eta_{D_{tors}}^\mathcal{E}$ is an isomorphism.
\endpf
\result{Remark}{\rm This proposition can be improved, see Remark~\ref{improved}.
}
\fresult
\vspace{.2cm}\par
Putting together the previous two propositions, we obtain the
following weak form of our main result.

\result{Theorem} \label{th1} Let $p$ be an odd prime number.
Let $\mathcal{X}$ be the class of $p$-groups consisting of all
elementary
abelian $p$-groups and 
the groups isomorphic to the extraspecial group of order $p^3$ and
exponent~$p$. Then the unit morphism
$$\eta_D^\mathcal{X} : D\to \mathcal{R}_\mathcal{X}\mathcal{O}_
\mathcal{X}D$$
is an isomorphism.
\fresult

\pf
We apply the functor $\mathcal{R}_\mathcal{X}\mathcal{O}_\mathcal{X}$
to the exact sequence
$$0\; \flh{}{} \, D_{tors} \, \flh{}{} D \flh{}{}  K^* \flh{}{} \; 0 $$
of Theorem~\ref{exact}. Since this functor is left exact
(because $\mathcal{O}_\mathcal{X}$ is exact and
$\mathcal{R}_\mathcal{X}$ has a left adjoint),
the commutative diagram
$$\xymatrix{
0\ar[r]&D_{tors}\ar[d]^{\eta_{D_{tors}}^\mathcal{X}}\ar[r]&D\ar[d]^
{\eta_D^\mathcal{X}}\ar[r]&K^*\ar[d]^{\eta_{K^*}^\mathcal{X}}\ar[r]&0\\
0\ar[r]&\mathcal{R}_\mathcal{X}\mathcal{O}_\mathcal{X}D_{tors}\ar[r]&
\mathcal{R}_\mathcal{X}\mathcal{O}_\mathcal{X}D\ar[r]&\mathcal{R}_
\mathcal{X}\mathcal{O}_\mathcal{X}K^*&
}$$
has exact rows.
By Propositions~\ref{K*iso} and~\ref{D-tors-iso}, both vertical arrows
$\eta_{D_{tors}}^\mathcal{X}$ and $\eta_{K^*}^\mathcal{X}$ are
isomorphisms.
Therefore $\eta_D^\mathcal{X}$ is an isomorphism.
\endpf

In order to pass from the class $\mathcal X$ to the smaller class
${\mathcal X}_3$, we need the following result about the Dade group of
elementary abelian $p$-groups.

\result{Lemma} \label{abel3} Let $\mathcal{E}_3$ be the class of all
elementary abelian $p$-groups of rank~$\leq3$.
If~$H$~is an elementary abelian $p$-group, the map
$$\eta^{\mathcal{E}_3}_D(H): D(H )\longrightarrow
\mathcal{R}_{\mathcal{E}_3}\mathcal{O}_{\mathcal{E}_3} D(H)$$
is an isomorphism.
\fresult

\pf
By Corollary~\ref{unit-defres}, the map $\eta^{\mathcal{E}_3}_D(H)$
coincides with the map
$$ \prod_{(T,S)\in\mathcal{E}_3(H)}\Defres_{T/S}^H : D(H )
\longrightarrow
\limproj{(T,S)\in\mathcal{E}_3(H)}D(T/S)\,.$$
We prove the result by induction on $|H|$. If $|H|\leq p^3$, then
$\limproj{(T,S)\in\mathcal{E}_3(H)}D(T/S) = D(H)$ and there is
nothing to prove. This starts induction. Assume $|H|>p^3$,
and let $u=(u_{T,S})_{(T,S)\in\mathcal{E}_3(H)}$ be an element
of $\limproj{(T,S)\in\mathcal{E}_3(H)}D(T/S)$.\par

Fix a non trivial subgroup $J$ of $H$. The sequence of elements
$u_{T,S}$, with $J\leq S\leq T\leq H$ and $|T/S|\leq p^3$ is an element
of $\limproj{(T,S)\in\mathcal{E}_3(H/J)}D(H/J)$.
Hence by induction  hypothesis, there exists a unique element
$v_{J}\in D(H/J)$ such that $\Defres_{T/S}^{H/J}v_J=u_{T,S}$,
for any section $(T,S)$ of $H$ with $J\leq S\leq T\leq H$
and $|T/S|\leq p^3$.\par

The uniqueness of $v_J$ implies that $\Def_{H/J'}^{H/J}v_J=v_{J'}$
whenever $J$ and $J'$ are subgroups of $H$ with $1<J\leq J'\leq H$.
So the sequence $(v_J)_{1<J\leq H}$ is an element of
$\limproj{1<J\leq H}D(H/J)$.
Now by Lemma~2.2 of~\cite{both2}, the deflation map
$$D(H) \longrightarrow \limproj{1<J\leq H}D(H/J)
$$
is an isomorphism. Explicitly, it was observed that the element
$$w=-\sum_{1<J\leq H}\mu(1,J)\Inf_{H/J}^Hv_J \,\in D(H)$$
is such that $\Def_{H/J}^Hw=v_J$ for any non trivial subgroup
$J$ of~$H$
(where $\mu$ denotes the M\"obius function of the poset of subgroups
of~$H$).
It follows that $\Defres_{T/S}^Hw=u_{T,S}$ for any
$(T,S)\in\mathcal{E}_3(H)$ with $S\neq1$.\par

Let $F$ be any subgroup of $H$ of order $p^3$, and consider
$f_F=\Res_F^Hw-u_{F,1}$.
Then clearly $\Def_{F/F'}^F f_F=0$ whenever $F'$ is a non trivial
subgroup of $F$.
It follows that
$$f_F\in\bigcap_{1<F'\leq F} \Ker(\Def_{F/F'}^F) \,,$$
which is the subgroup of $D(F)$ consisting of endo-trivial modules
(see Lemma~1.2 in~\cite{both1}).
Since this subgroup is infinite cyclic generated by the
class~$\Omega_F$ of the augmentation ideal of~$kF$
(Dade's theorem, see Theorem~1.4 in~\cite{both1}), there exists a
unique integer $m_F$ such that $f_F=m_F\Omega_{F}$.
Now if $F'$ is another subgroup of order $p^3$ of $H$
such that $F\cap F'$ has order at least $p^2$, we have

\begin{eqnarray*}
\Res_{F\cap F'}^Ff_F&=&\Res_{F\cap F'}^Hw-\Res_{F\cap F'}^Fu_{F,1}\\
&=&\Res_{F\cap F'}^Hw-u_{F\cap F',1}\\
&=&\Res_{F\cap F'}^{F'}f_{F'}\,.
\end{eqnarray*}
Thus $m_F\Omega_{F\cap F'}=m_{F'}\Omega_{F\cap F'}$, and $m_F=m_{F'}$
because $\Omega_{F\cap F'}$ has infinite order.
Thus $m_F=m_{F'}$ if $|F\cap F'|\geq p^2$, and since the poset of
subgroups of $H$ of order $p^2$ and~$p^3$ is connected, $m_F$ does
not depend on the subgroup $F$ of order $p^3$. Set $m=m_F$, and
consider the element $t=w-m\Omega_{H}$. For any section $(T,S)\in
\mathcal{E}_3(H)$ with $S\neq 1$, we have again
$$\Defres_{T/S}^Ht=\Defres_{T/S}^Hw=u_{T,S}\,.$$
Moreover, if $E$ is any subgroup of $H$ of order at most $p^3$,
choose some subgroup $F$ of order $p^3$ containing~$E$. Then
$$\Defres_{E/1}^Ht=\Res_E^Ht=\Res_E^F\Res_F^Ht=\Res_E^Fu_{F,1}=u_{E,1}
\mpoint$$
It follows that $\eta^{\mathcal{E}_3}_D(H)(t)=u$, and so
$\eta^{\mathcal{E}_3}_D(H)$ is surjective.
It is also injective by Theorem~\ref{detection} and
this completes the proof of the lemma.
\endpf

We have now paved the way for the final version of our main result.

\result{Theorem} \label{th2}
Let $p$ be an odd prime number and let $\mathcal{X}_3$ be the class
of $p$-groups consisting of all elementary abelian $p$-groups
of rank~$\leq3$ and 
the groups isomorphic to the extraspecial group of order $p^3$
and exponent~$p$.
Then the unit morphism
$$\eta_D^{\mathcal{X}_3} : D\to \mathcal{R}_{\mathcal{X}_3}\mathcal{O}
_{\mathcal{X}_3}D$$
is an isomorphism. In other words, if $P$ is a $p$-group, then the map
$$\eta_{D}^{\mathcal{X}_3}(P) \;=\! \prod_{(T,S)\in{\mathcal{X}_3}(P)}
\Defres_{T/S}^P : \; D(P)\longrightarrow \limproj{(T,S)\in{\mathcal{X}
_3}(P)} D(T/S)$$
is a group isomorphism.
\fresult

\pf
We are going to apply Lemma~\ref{rapetisser} to the class~$\mathcal X
$ and
the subclass ${\mathcal{X}_3}$. First note that the unit morphism $
\eta_D^{\mathcal{X}}$ is an isomorphism by Theorem~\ref{th1}.
Moreover, we must show that, for any $Q\in\mathcal{X}$, the unit
morphism $\eta_D^{\mathcal{X}_3}(Q)$ evaluated at~$Q$ is an isomorphism.
This is obvious if $Q\in\mathcal{X}_3$ (in that case $\eta_D^{\mathcal
{X}_3}(Q)$ is the identity).
Now if $Q\in\mathcal{X}-\mathcal{X}_3$, then $Q$ is elementary
abelian of rank~$\geq4$ and the family $\mathcal{X}_3(Q)$ coincides
with $\mathcal{E}_3(Q)$.
Thus $\eta_D^{\mathcal{X}_3}(Q)=\eta_D^{\mathcal{E}_3}(Q)$ and this
is an isomorphism by Lemma~\ref{abel3}.
Thus the assumptions of Lemma~\ref{rapetisser} are satisfied and
therefore $\eta_D^{\mathcal{X}_3}$ is an isomorphism.
\endpf

\result{Remark} \label{improved}
 {\rm Proposition~\ref{D-tors-iso} can be improved in the same way, replacing the class~$\mathcal{E}$ by the smaller class $\mathcal{E}_3$. Thus we obtain that,
if $p$ is odd, the unit homomorphism
$$\eta_{D_{tors}}^{\mathcal{E}_3} : D_{tors}\longrightarrow
\mathcal{R}_{\mathcal{E}_3}\mathcal{O}_{\mathcal{E}_3}D_{tors}$$
is an isomorphism.
}
\fresult 

\result{Remark} {\rm The class $\mathcal{X}_3$ is exactly the class
of $p$-groups of order at most $p^3$ and exponent dividing~$p$.
}
\fresult

\result{Remark} {\rm Theorem~\ref{th2} holds more generally for any
class of $p$-groups containing~$\mathcal{X}_3$.
This follows immediately from Lemma~\ref{agrandir}.}
\fresult

\result{Remark} {\rm Theorem~\ref{th2} does not hold if we remove
from $\mathcal{X}_3$ 
the groups which are isomorphic to one of the two groups of order~$p^3$ (elementary
abelian or extraspecial). In order to prove this,
we let $P\in\mathcal{X}_3$ with $|P|=p^3$,
we consider the unit morphism $\eta_D^{\mathcal{E}_2}(P)$
corresponding to proper sections of~$P$ and we show that
$\eta_D^{\mathcal{E}_2}(P)$ is not surjective. Note first that $P$ has
more than one elementary abelian subgroup of rank~2, and that each of
them is normal in $P$.
Let $Q$ be any of them and define, for any $(T,S)\in\mathcal{E}_2(P)$,
$$u_{T,S}=\left\{\begin{array}{ll}2\Omega_Q&\hbox{if}\;(T,S)=(Q,1)\;,\\
0&\hbox{otherwise}\;.\end{array}\right.$$
The family $u=(u_{T,S})$ defines an element of
$\limproj{(T,S)\in{\mathcal{E}_2}(P)} D(T/S)$ because
$\Res_R^Q(2\Omega_Q)=2\Omega_R=0$ whenever $R<Q$
(because $D(R)=\Z/2\Z$ if $|R|=p$).
We want to prove that $u$ is not in the image
of $\eta_D^{\mathcal{E}_2}(P)$.
Suppose there exists $v\in D(P)$ whose image is~$u$. Then any section
$(T,S)$ of~$P$ with $S\neq 1$ belongs to ${\mathcal{E}_2}(P)$ and
$\Defres_{T/S}^P(v)=u_{T,S}=0$ by definition of~$u$.
It follows that $v$ belongs to the subgroup~$T(P)$ of~$D(P)$ consisting
of classes of endo-trivial modules.
Moreover $\Res_Q^P(v)=2\Omega_Q$ and $\Res_{Q'}^P(v)=\nolinebreak 0$
for any other
elementary abelian subgroup~$Q'$ of rank~2.
This is impossible if $P$ is elementary abelian of rank~3 because
$T(P)=\Z$
generated by~$\Omega_P$ and we would obtain
$v=2\Omega_P$ (because $\Res_Q^P(v)=2\Omega_Q$)
and also $v=0$ (because $\Res_{Q'}^P(v)=0$). This is also impossible
if $P$ is extraspecial of order~$p^3$ and exponent~$p$ because the
condition $\Res_{Q'}^P(v)=0$ for any $Q'\neq Q$ implies that
$\Res_Q^P(v)$ must be a multiple of $2p\Omega_Q$, using the
description of $T(P)$ given in~\cite{cath3} (see Theorems~3.1 and 6.1
of that paper).
}
\fresult

\section{Appendix}

\noindent
The purpose of this appendix is to provide a proof of
Theorem~\ref{adjoints}.
We shall only prove the statement about $\mathcal{R}_\mathcal{Y}$
because the treatment of $\mathcal{L}_\mathcal{Y}$ is
similar.
This type of construction of adjoints is more or less standard in
category theory (see for instance~\cite{GZ}).
We provide here an explicit treatment in our situation.\par

We first start with technical lemmas about bisets, using the
following notation.

\result{Notation} {\rm Let $P$ and $Q$ be two groups.
\begin{itemize}
\item[(a)] If $V$ is a right $Q$-set and $U$ is a left $Q$-set, $(v,_
{_Q}u)$ denotes the image
in $V\times_QU$ of the pair $(v,u)$ of $V\times U$.
\item[(b)] Let $U$ be a $(Q,P)$-biset. If $S$ is a subgroup of $P$
and $u$ is an element of $U$, set
$$\ls uS=\{y\in Q\mid\exists s\in S,\;us=yu\}\mpoint$$
Then $\ls uS$ is a subgroup of $Q$. Similarly, if $T$ is a subgroup
of $Q$, set
$$T^u=\{x\in P\mid\exists t\in T,\;tu=ux\}\mpoint$$
It is a subgroup of $P$.
\end{itemize}
}
\fresult

\result{Remark} {\rm
The notation in~(b) extends the standard notation for conjugation in
the following sense. Let $P$ be a $p$-group and $S$ a subgroup of~$P$.
\begin{itemize}
\item[(1)] If $P$ is a subgroup of~$Q$ and $U=P$ viewed as
a $(Q,P)$-biset,
then $\ls uS = uSu^{-1}$.
\item[(2)] If $Q$ is a subgroup of~$P$ and $U=P$ viewed as
a $(Q,P)$-biset,
then $\ls uS= Q\cap uSu^{-1}$.
\item[(3)] If $Q=P/R$ is a quotient group of~$P$ and $U=P/R$ viewed
as a $(Q,P)$-biset,
then $\ls uS$ is the image in~$Q$ of the conjugate subgroup~$uSu^{-1}$.
\item[(4)] If $P=Q/R$ is a quotient group of~$Q$ and $U=Q/R$ viewed
as a $(Q,P)$-biset,
then $\ls uS$ is the inverse image in~$Q$ of the conjugate
subgroup~$uSu^{-1}$.
\end{itemize}
}
\fresult

\result{Remark} {\rm Note that $1^u$ is the stabilizer of~$u$ in~$P$
and $\ls u1$ is the stabilizer of~$u$ in~$Q$.
In the definition of~$\ls uS$, the element $s$ is unique up to left
multiplication by an element of~$S\cap 1^u$.
Similarly, in the definition of~$T^u$, the
element $t$ is unique up to right multiplication by an element
of~$T\cap \ls u1$.
This explains why these subgroups appear in parts~(b) and (b') of the
next lemma.

}
\fresult

\result{Lemma} \label{conjugaison} Let $P$, $Q$ and $R$ be groups,
let $U$ be a $(Q,P)$-biset, and let $V$ be an $(R,Q)$-biset.
Let $u\in U$ and $v\in V$.
\begin{itemize}
\item[(a)] If $T$ is a subgroup of $Q$ and if $x\in P$,
then $(T^u)^x=T^{ux}$.
\item[(a')] If $X$ is a subgroup of $P$ and if $y\in Q$,
then $^y({\ls uX})={\ls {yu}X}$.
\item[(b)] If $(T,S)$ is a section of $Q$, then $(T^u,S^u)$ is a
section of $P$, and
there are group isomorphisms
$$T^u/S^u\cong (T\cap {\ls uP})/(S\cap {\ls uP})(T\cap{\ls u1})\cong
(T\cap{\ls uP})S/(T\cap{\ls u1})S\mpoint$$
In particular, the quotient $T^u/S^u$ is isomorphic to
a subquotient of~$T/S$.
\item[(b')] If $(Y,X)$ is a section of $P$, then $({\ls uY},{\ls uX})
$ is a section of $Q$, and there are group isomorphisms
$${\ls uY}/{\ls uX}\cong (Y\cap Q^u)/(X\cap Q^u)(Y\cap 1^u)\cong (Y
\cap Q^u)X/(Y\cap 1^u)X\mpoint$$
In particular, the quotient ${\ls uY}/{\ls uX}$ is isomorphic to a
subquotient of $Y/X$.
\item[(c)] If $X$ is a subgroup of $P$, then $\ls v({\ls uX})={\ls
{(v,_{_Q}u)}X}$.
\item[(c')] If $Z$ is a subgroup of $R$, then $(Z^v)^u=Z^{(v,_{_Q}u)}$.
\end{itemize}
\fresult

\pf We only prove (a), (b), (c), because the proofs of (a'), (b'),
(c') are similar.\par

(a) By definition,
\begin{eqnarray*}
T^{ux}&=&\{g\in P\mid\exists t\in T,\;tux=uxg\}\\
&=&\{g\in P\mid\exists t\in T,\;tu=uxgx^{-1}\}\\
&=&\{g\in P\mid {\ls xg}\in T^u\}=(T^u)^x\mpoint
\end{eqnarray*}
\par

(b) Suppose that $(T,S)$ is a section of $Q$. Then $S^u$ is clearly a
subgroup of~$T^u$.
Moreover if $x\in T^u$, then there exists $t\in T$ such that $tu=ux$.
Thus
$$(S^u)^x=S^{ux}=S^{tu}=S^u\mpoint$$
Here the first equality follows from (a) and the last one holds
because $t\in T$ normalizes~$S$. This shows that $S^u\normal T^u$.
Consider the map
$$T^u/S^u\longrightarrow (T\cap{\ls uP})/(S\cap{\ls uP})(T\cap{\ls
u1}) \,, \qquad
xS^u \mapsto t(S\cap{\ls uP})(T\cap{\ls u1}) \,,$$
where $t$ is any element of $T$ such that $tu=ux$.
It is routine to check that it is a well defined group homomorphism.
Similarly the map
$$(T\cap{\ls uP})/(S\cap{\ls uP})(T\cap{\ls u1}) \longrightarrow T^u/
S^u \,, \qquad
t(S\cap{\ls uP})(T\cap{\ls u1}) \mapsto xS^u \,,$$
where $x$ is any element of $P$ such that $tu=ux$, is a well defined
group homomorphism. Clearly these two group homomorphisms are mutual
inverse.\par

Now setting $N=(T\cap{\ls u1})S$, we have
$$(T\cap{\ls uP})S/(T\cap{\ls u1})S=(T\cap{\ls uP})N/N\cong (T\cap
{\ls uP})/(T\cap{\ls uP}\cap N)\mvirg$$
and
$$T\cap{\ls uP}\cap N=T\cap{\ls uP}\cap (T\cap{\ls u1})S=(T\cap{\ls
u1})(T\cap{\ls uP}\cap S)=(T\cap{\ls u1})(S\cap {\ls uP})\mpoint$$
This proves the second isomorphism.\par

Now the group $(T\cap{\ls uP})/(S\cap{\ls uP})(T\cap{\ls u1})$ is a
factor group of $(T\cap{\ls uP})/(S\cap{\ls uP})$, which is
isomorphic to $(T\cap{\ls uP})S/S$, and this is a subgroup of $T/S$.
Thus $T^u/S^u$ is isomorphic to a subquotient of $T/S$.
\par

(c) By definition
\begin{eqnarray*}
^{(v,_{_Q}u)}X&=&\{z\in R\mid\exists x\in X,\;z(v,_{_Q}u)=(v,_{_Q}u)x
\}\\
&=&\{z\in R\mid\exists x\in X,\;(zv,_{_Q}u)=(v,_{_Q}ux)\}\\
&=&\{z\in R\mid\exists x\in X,\;\exists y\in Q,\;zv=vy,\;yu=ux\}\\
&=&\{z\in R\mid\exists y\in Q,\;\exists x\in X,\;zv=vy,\;yu=ux\}\\
&=&\{z\in R\mid\exists y\in {\ls uX},\;zv=vy\}\\
&=&{\ls v({\ls uX})}
\end{eqnarray*}
This completes the proof.
\endpf

\result{Lemma} \label{facile}
Let $P$, $Q$, and $R$ be groups, let $U$ be a $(Q,P)$-biset, and let
$V$ be an $(R,Q)$-biset. Let moreover $(D,C)$ be a section of~$R$,
let $(B,A)$ be a section of~$Q$ and
assume that $A$ acts trivially on~$C\bs V$ (on the right).
Then there is an isomorphism of $(D/C,P)$-bisets
$$(C\bs V)\times_{B}(A\bs U)\cong C\bs (V\times_BU)\mpoint$$
\fresult

\pf The maps $(Cv,_{_B}Au)\mapsto C(v,_{_B}u)$ and $C(v,_{_B}u)
\mapsto (Cv,_{_B}Au)$ are well defined, and mutual inverse biset
isomorphisms.
\endpf

\bigskip

Let $\mathcal{Y}$ be a class of finite $p$-groups, closed under
taking sections. If $F$ is biset functor defined on~$\mathcal{Y}$ and
if $P$ is an arbitrary $p$-group, we have defined

\begin{equation}\label{lim}
\mathcal{R}_{\mathcal{Y}}F(P)=\limproj{(T,S)\in \mathcal{Y}(P)}F(T/S)
\mpoint
\end{equation}

Recall that $\mathcal{Y}(P)$ is the set of sections $(T,S)$ of $P$
such that $T/S\in\mathcal{Y}$ and that the definition of the limit
means that the group $\mathcal{R}_{\mathcal{Y}}F(P)$ is the set of
sequences
$(l_{T,S})_{(T,S)\in\mathcal{Y}(P)}$ indexed by $\mathcal{Y}(P)$,
where $l_{T,S}\in F(T/S)$, subject to the following conditions~:
\begin{enumerate}
\item If $(T,S)$ and $(T',S')$ are elements of $\mathcal{Y}(P)$ such
that $S\leq S'\leq T'\leq T$, then
$$\Defres_{T'/S'}^{T/S}l_{T,S}=l_{T',S'}\mvirg$$
where $\Defres_{T'/S'}^{T/S}$ is the set $S'\bs T$, viewed as
a $(T'/S',T/S)$-biset.
\item If $(T,S)\in \mathcal{Y}(P)$ and if $x\in P$, then
$$^xl_{T,S}=l_{\ls xT,{\ls xS}}\mpoint$$
\end{enumerate}
Now we want to introduce a structure of biset functor
on $\mathcal{R}_{\mathcal{Y}}F$.
Whenever $P$ and $Q$ are $p$-groups and $U$ is a finite
$(Q,P)$-biset, we need to define a map
$$\mathcal{R}_{\mathcal{Y}}F(U): \mathcal{R}_{\mathcal{Y}}F(P)\to
\mathcal{R}_{\mathcal{Y}}F(Q) \,.$$
If $l$ is an element of $\mathcal{R}_{\mathcal{Y}}F(P)$, and if
$(T,S)\in\mathcal{Y}(Q)$, we set

\begin{equation}\label{rfu}
\mathcal{R}_{\mathcal{Y}}F(U)(l)_{T,S}=\sum_{u\in [T\bs U/P]}(S\bs Tu)
\cdot l_{T^u,S^u}\mvirg
\end{equation}

\noindent
where $[T\bs U/P]$ is any set of representatives of $(T,P)$-orbits on
$U$, where $S\bs Tu$ is viewed as a $(T/S,T^u/S^u)$-biset, and
$(S\bs Tu)\cdot l_{T^u,S^u}$ denotes the image of $l_{T^u,S^u}$
by this biset. This makes sense since $(T^u,S^u)\in\mathcal{Y}(P)$
if $(T,S)\in\mathcal{Y}(Q)$, by Lemma~\ref{conjugaison}.

\result{Remark} {\rm If we want to make this definition explicit in
terms of elementary operations, we observe that it is equivalent to}
\begin{equation*}
\mathcal{R}_{\mathcal{Y}}F(U)(l)_{T,S}=\sum_{u\in [T\bs U/P]}
\Indinf_{(T\cap {\ls uP})S/(T\cap{\ls u1})S}^{T/S}\Iso^{(T\cap{\ls
uP})S/(T\cap{\ls u1})S}_{T^u/S^u}l_{T^u,S^u}\mpoint
\end{equation*}
\fresult

\result{Lemma} The element $\mathcal{R}_{\mathcal{Y}}F(U)(l)_{T,S}$
is well defined, that is, it does not depend on the choice of the set
of representatives $[T\bs U/P]$.
\fresult

\pf Changing the set of representatives $[T\bs U/P]$ amounts to
replacing each element~$u$ in this set by an element $u'=t_uux_u$,
where $t_u\in T$ and $x_u\in P$. This gives a new set of
representatives of orbits, denoted by $[T\bs U/P]'$. With this new
set of representatives, the sum in~(\ref{rfu}) becomes
\begin{eqnarray*}
\Sigma{'}&=&\sum_{u'\in [T\bs U/P]'}(S\bs Tu')\cdot l_{T^{u'},S^{u'}}\\
&=&\sum_{u\in [T\bs U/P]}(S\bs Tux_u)\cdot l_{(T^{u})^{x_u},(S^{u})^
{x_u}}\\
&=&\sum_{u\in [T\bs U/P]}(S\bs Tux_u)\cdot \ls{x_u^{-1}} l_{T^u,S^u}\;
\;\;
{\rm (since }\;\;\;l\in \mathcal{R}_{\mathcal{Y}}F(P){\rm ),}\\
&=&\sum_{u\in [T\bs U/P]}(S\bs Tu)\cdot l_{T^u,S^u}\mpoint
\end{eqnarray*}
The latter equality follows from the fact that, for any $x\in P$, $(S
\bs Tu) \circ \Conj(x) \cong S\bs Tux$ as
$(T/S,T^{ux}/S^{ux})$-bisets
(where $\Conj(x)$ is the $(T^u/S^u,T^{ux}/S^{ux})$-biset
$T^u/S^u$ with right action of $T^{ux}/S^{ux}$ consisting of conjugation
by~$x$ followed by right multiplication).
The lemma follows.
\endpf

\result{Example} \label{ex1}
{\rm Let $F\in\mathcal{F}_{\mathcal{Y}}$, let $P$ be a $p$-group, let
$(T,S)$ and $(T',S')$ be elements of $\mathcal{Y}(P)$ such that
$S\leq S'\leq T'\leq T$. We view $(T',S')$ as a section of~$T/S$ via
the canonical isomorphism $(T'/S)/(S'/S) \cong T'/S'$.
If $l\in\mathcal{R}_{\mathcal{Y}}F(P)$, then
$$(\Defres_{T/S}^Pl)_{T',S'}=l_{T',S'}\mpoint$$
Indeed in this case, the group $Q$ is equal to $T/S$ and the biset
$U$ is~$S\bs P$. Hence $[T'\bs U/P]$ has one element, which can be
chosen to be $\{S\}$, and in this case
$$\mathcal{R}_{\mathcal{Y}}F(U)(l)_{T',S'}=(S'\bs T')\cdot l_{T',S'}
=l_{T',S'}\mvirg$$
because the $(T'/S',T'/S')$-biset $S'\bs T'$ is the identity.
}
\fresult

\result{Lemma} \label{family}
Let $l\in \mathcal{R}_{\mathcal{Y}}F(P)$.
When $(T,S)$ runs through~${\mathcal{Y}}(Q)$,
the family $\mathcal{R}_{\mathcal{Y}}F(U)(l)_{T,S}$
defined in~(\ref{rfu}) is an element
of~$\mathcal{R}_{\mathcal{Y}}F(Q)$.
\fresult

\pf Let $(T,S)$ and $(T',S')$ be elements of $\mathcal{Y}(Q)$ such
that $S\leq S'\leq T'\leq T$. Then $\Defres_{T'/S'}^{T/S}$ is the set
$S'\bs T$, viewed as a $(T'/S',T/S)$-biset. It follows that
\begin{eqnarray*}
\Defres_{T'/S'}^{T/S}\bigl(\mathcal{R}_{\mathcal{Y}}F(U)(l)_{T,S}\bigr)
&=&S'\bs T\times_{T/S}\sum_{u\in [T\bs U/P]}(S\bs Tu)\cdot l_{T^u,S^u}\\
&=&\sum_{u\in [T\bs U/P]}(S'\bs Tu)\cdot l_{T^u,S^u}\mvirg\end
{eqnarray*}
because by Lemma~\ref{facile}, there is a biset isomorphism
$$(S'\bs T)\times_{T/S}(S\bs Tu)\cong S'\bs (T\times_TTu)
=S'\bs Tu\mpoint$$
Now the set of orbits of $(T'/S')\times (T^u/S^u)$ on the biset
$S'\bs Tu$ is the set $T'\bs Tu/T^u$, which is in bijective
correspondence with the set $T'\bs T/T\cap {\ls uP}$. Hence there is
an isomorphism of $(T'/S',T^u/S^u)$-bisets
$$S'\bs Tu\cong \bigsqcup_{t\in[T'\bs T/T\cap {\ls uP}]}S'\bs T'tuT^u
\mpoint$$
It follows that

\begin{equation}\label{one-side}
\Defres_{T'/S'}^{T/S}\bigl(\mathcal{R}_{\mathcal{Y}}F(U)(l)_{T,S}\bigr)
=\sumb{u\in [T\bs U/P]}{t\in[T'\bs T/T\cap {\ls uP}]}S'\bs T'tuT^u
\cdot l_{T^u,S^u}\mpoint
\end{equation}
On the other hand, we have by definition
\begin{eqnarray*}
\mathcal{R}_{\mathcal{Y}}F(U)(l)_{T',S'}&=&\sum_{v\in [T'\bs U/P]}
(S'\bs T'v)\cdot l_{T'^v,S'^v}\\
&=&\sumb{u\in [T\bs U/P]}{t\in [T'\bs T/T\cap {\ls uP}]}(S'\bs T'tu)
\cdot l_{T'^{tu},S'^{tu}}
\end{eqnarray*}
Now $S^u\leq S'^{tu}\leq T'^{tu}\leq T^u$, thus by Example~\ref{ex1}
$$l_{T'^{tu},S'^{tu}}=\Defres_{T'^{tu}/S'^{tu}}^{T^u/S^u}\;l_{T^u,S^u}
=S'^{tu}\bs T^u\cdot l_{T^u,S^u}\mpoint$$
It follows that
\begin{eqnarray*}
\mathcal{R}_{\mathcal{Y}}F(U)(l)_{T',S'}
&=&\!\!\!\sumb{u\in [T\bs U/P]}{t\in [T'\bs T/T\cap {\ls uP}]}
\!\!\!(S'\bs T'tu)\times_{T'^{tu}/S'^{tu}}(S'^{tu}\bs T^u)\cdot l_
{T^u,S^u}\mpoint\\
\end{eqnarray*}
We have to prove that this sum coincides with the expression
in~(\ref{one-side}). By Lemma~\ref{facile},
for any $u\in U$ and any $t\in T$,
there is an isomorphism of $(T'/S',T^u/S^u)$-bisets
$$(S'\bs T'tu)\times_{T'^{tu}/S'^{tu}}(S'^{tu}\bs T^u)
\cong S'\bs (T'tu\times_{T'^{tu}}T^u) \mpoint$$
But it is easy to see that
$T'tu\times_{T'^{tu}}T^u \cong T'tuT^u$ as $(T',T^u)$-bisets, so that
we obtain
$S'\bs (T'tu\times_{T'^{tu}}T^u) \cong S'\bs T'tu(T^u)$.
\endpf

\result{Lemma} \label{composition}
Let $P$, $Q$, and $R$ be finite $p$-groups.
\begin{itemize}
\item[(a)] Let $U$ be the set $P$, viewed as a $(P,P)$-biset by left
and right multiplication.
Then $\mathcal{R}_{\mathcal{Y}}F(U)$ is the identity map.
\item[(b)]  Let $U$ be a $(Q,P)$-biset and let $V$ be
an $(R,Q)$-biset. Then
$$\mathcal{R}_{\mathcal{Y}}F(V)\circ \mathcal{R}_{\mathcal{Y}}F(U)=
\mathcal{R}_{\mathcal{Y}}F(V\times_QU)\mpoint$$
\end{itemize}
\fresult

\pf (a) Let $l\in \mathcal{R}_{\mathcal{Y}}F(P)$, and $(T,S)\in
\mathcal{Y}(P)$. The set $[T\bs U/P]$ can be chosen to be equal
to $\{1\}$, and for $u=1$, one has that $T^u=T$ and $S^u=S$.
Thus
$$\mathcal{R}_{\mathcal{Y}}F(U)(l)_{T,S}=(S\bs T)\cdot l_{T,S}=l_{T,S}
\mvirg$$
since $S\bs T$ acts as the identity on $F(T/S)$.\par
(b) Let $l\in \mathcal{R}_{\mathcal{Y}}F(P)$ and $(T,S)\in \mathcal{Y}
(R)$, and denote by $L$ the element $\big(\mathcal{R}_{\mathcal{Y}}F
(V)\circ \mathcal{R}_{\mathcal{Y}}F(U)(l)\big)_{T,S}$. Then
\begin{eqnarray*}
L&=&\sum_{v\in [T\bs V/Q]}(S\bs Tv)\cdot \mathcal{R}_{\mathcal{Y}}F(U)
(l)_{T^v,S^v}\\
&=&
\sumb{v\in [T\bs V/Q]}{u\in [T^v\bs U/P]}(S\bs Tv)\times_{T^v/S^v}(S^v
\bs T^vu)\cdot l_{(T^v)^u,(S^v)^u}
\end{eqnarray*}
But by Lemma~\ref{conjugaison}, $(T^v)^u=T^{(v,_{_Q}u)}$ and
$(S^v)^u=S^{(v,_{_Q}u)}$.
Moreover when $v$ runs through $[T\bs V/Q]$ and
$u$ runs through $[T^v\bs U/P]$, the element $(v,_{_Q}u)$
runs through
a set of representatives of the set $R\bs (V\times_QU)/P$, and by
Lemma~\ref{facile}, there is an isomorphism
of $(T/S,T^{(v,_{_Q}u)})$-bisets
$$(S\bs Tv)\times_{T^v/S^v}(S^v\bs T^vu)\cong S\bs (Tv\times_{T^v}
T^vu)\mpoint$$
Moreover, there is an isomorphism of $(T,T^{(v,_{_Q}u)})$-bisets
$$Tv\times_{T^v}T^vu\cong T(v,_{_{Q}}u) \,,$$
sending $(tv,_{_{T^v}}xu)$ to $tt'(v,_{_Q}u)$, for $t\in T$ and
$x\in T^v$, where $t'\in T$ is such that $t'v=vx$.
The inverse isomorphism
maps $t(v,_{_{Q}}u)$ to $(tv,_{_{T^v}}u)$, for $t\in T$.
In order to see that this is well-defined, replace $t$ by $tw$ where
$w\in T$ stabilizes $(v,_{_Q}u)$.
Then $(wv,_{_Q}u)=(v,_{_Q}u)$, so that there exists $y\in Q$ such
that $wv=vy$ and $u=y^{-1}u$.
We now see that $y\in T^v$ and it follows that
$$(twv,_{_{T^v}}u)=(tvy,_{_{T^v}}u)=(tv,_{_{T^v}}yu)=(tv,_{_{T^v}}u)
\mpoint$$
Thus finally
\begin{eqnarray*}
L&=&\sum_{(v,_{_Q}u)\in[R\bs (V\times_QU)/P]}\big(S\bs T(v,_{_Q}u)
\big)\cdot l_{T^{(v,_{_Q}u)},S^{(v,_{_Q}u)}}\\
&=&\mathcal{R}_{\mathcal{Y}}F(V\times_QU)(l)_{T,S}\mvirg
\end{eqnarray*}
and the result follows.
\endpf

\result{Theorem} Let $F$ be a biset functor defined
on~${\mathcal{Y}}$.
The correspondences mapping a $p$-group~$P$
to $\mathcal{R}_{\mathcal{Y}}F(P)$ and
a $(Q,P)$-biset~$U$ to
$\mathcal{R}_{\mathcal{Y}}F(U)$ define a functor
from $\mathcal{C}_p$ to $\mathcal{A}b$ (in
other words a biset functor on $p$-groups).
\fresult

\pf This follows from Lemmas~\ref{family} and \ref{composition},
together with the observation that the correspondence $U\mapsto
\mathcal{R}_{\mathcal{Y}}F(U)$ is obviously additive with respect to
the biset~$U$, in the sense that if $U$ is a disjoint union~$U=U_1
\sqcup U_2$, then
$\mathcal{R}_{\mathcal{Y}}F(U)=
\mathcal{R}_{\mathcal{Y}}F(U_1) + \mathcal{R}_{\mathcal{Y}}F(U_2)$.
\endpf

Recall that $\mathcal{F}$ denotes the category of biset functors
(defined on all $p$-groups)
and $\mathcal{F}_{\mathcal{Y}}$ the category of biset functors
defined on~$\mathcal{Y}$.
Moreover $\mathcal{O} : \mathcal{F} \to \mathcal{F}_{\mathcal{Y}}$
denotes the forgetful functor.

\result{Theorem} \label{right-adj}
The correspondence $F\mapsto \mathcal{R}_{\mathcal{Y}}F$ is a functor
$\mathcal{R}_{\mathcal{Y}}$ from $\mathcal{F}_{\mathcal{Y}}$
to $\mathcal{F}$, and this functor is right adjoint to the forgetful
functor $\mathcal{O}_{\mathcal{Y}}$.
Moreover the composition
$\mathcal{O}_{\mathcal{Y}}\circ\mathcal{R}_{\mathcal{Y}}$ is
isomorphic to the identity functor of $\mathcal{F}_{\mathcal{Y}}$.
\fresult

\pf We first show that the correspondence
$F\mapsto \mathcal{R}_{\mathcal{Y}}F$ is functorial in~$F$.
A morphism $\varphi: F\to G$ in $\mathcal{F}_{\mathcal{Y}}$ is a
collection of morphisms $\varphi_Q: F(Q)\to G(Q)$,
for $Q$ in~$\mathcal{Y}$, subject to the usual naturality conditions.
(Note that throughout this proof, we write $\varphi_Q$ instead of
$\varphi(Q)$ for simplicity of notation.)
If $P$
is a $p$-group and $l\in \mathcal{R}_{\mathcal{Y}}F(P)$,
and if $(T,S)\in\mathcal{Y}(P)$, we set
$$l'_{T,S}=\varphi_{T/S}(l_{T,S})\mpoint$$
Since the maps $\varphi_{T/S}$ commute with deflation maps,
restriction maps, and conjugation maps by the above naturality
conditions, it is it clear that this defines an element
$l'\in\mathcal{R}_{\mathcal{Y}}G(P)$.
Now the correspondence $l\mapsto l'$ is a linear map
$$\mathcal{R}_{\mathcal{Y}}(\varphi)_P:\;\mathcal{R}_{\mathcal{Y}}F(P)
\to \mathcal{R}_{\mathcal{Y}}G(P)\mpoint$$
Suppose now that $P$ and $Q$ are $p$-groups, and that $U$ is
a finite $(Q,P)$-biset. With the previous notation,
\begin{eqnarray*}
\mathcal{R}_{\mathcal{Y}}G(U)(l')_{T,S}&=&\sum_{u\in[T\bs U/P]}(S\bs
Tu)\cdot \varphi_{T^u/S^u}(l_{T^u,S^u})\\
&=&\sum_{u\in[T\bs U/P]}\varphi_{T/S}\big((S\bs Tu)\cdot l_{T^u,S^u}
\big)\\
&=&\varphi_{T/S}\big(\mathcal{R}_{\mathcal{Y}}F(U)(l)_{T,S}\big)
\end{eqnarray*}
This shows that $\mathcal{R}_{\mathcal{Y}}G(U)\circ \mathcal{R}_
{\mathcal{Y}}(\varphi)_P=\mathcal{R}_{\mathcal{Y}}(\varphi)_Q\circ
\mathcal{R}_{\mathcal{Y}}F(U)$.
In other words, the maps
$\mathcal{R}_{\mathcal{Y}}(\varphi)_P$
define a natural transformation of functors
$\mathcal{R}_{\mathcal{Y}}(\varphi):\;\mathcal{R}_{\mathcal{Y}}F\to
\mathcal{R}_{\mathcal{Y}}G$. It is now clear that the correspondence
$F\mapsto \mathcal{R}_{\mathcal{Y}}F$ is a functor
$\mathcal{F}_{\mathcal{Y}}\to \mathcal{F}$.\par

It is convenient to prove the last assertion of the theorem now. We
have to show that if $F\in\mathcal{F}_{\mathcal{Y}}$ and if $Q$ is in
the class~$\mathcal{Y}$, then there is an isomorphism
$$\varepsilon_Q:\;F(Q)\to \mathcal{O}_{\mathcal{Y}} \mathcal{R}_
{\mathcal{Y}}F(Q)$$
which is natural with respect to $F$ and~$Q$. This isomorphism is
defined as follows. If $f\in F(Q)$ and $(T,S)\in\mathcal{Y}(Q)$, then
$\varepsilon_Q(f)_{T,S}=\Defres_{T/S}^Qf$. It is routine to check
that this definition gives the required isomorphism of functors
$\Id\to \mathcal{O}_{\mathcal{Y}}\circ\mathcal{R}_{\mathcal{Y}}$.\par

We turn now to the adjointness property.
Suppose that $F\in \mathcal{F}$, $G\in\mathcal{F}_{\mathcal{Y}}$,
and $\varphi:~\mathcal{O}_{\mathcal{Y}}F\to G$ is a morphism
in~$\mathcal{F}_{\mathcal{Y}}$.
Thus $\varphi$ is a collection of morphisms $\varphi_Q:F(Q)\to G(Q)$,
for $Q\in\mathcal{Y}$, satisfying some commutation conditions.
Now if $P$ is any $p$-group, $f\in F(P)$, and
$(T,S)\in \mathcal{Y}(P)$, set
$$l_{T,S}=\varphi_{T/S}(\Defres_{T/S}^Pf)\,.$$
The above commutation conditions show easily that these elements
$l_{T,S}$ define
an element $l$ of~$\mathcal{R}_{\mathcal{Y}}G(P)$.
Now the correspondence $f\mapsto l$ is a linear map
from $F(P)$ to $\mathcal{R}_{\mathcal{Y}}G(P)$,
denoted by $\varphi^+_P$, such that
\begin{equation} \label{adjoint-defres}
\varphi^+_P(f)_{T,S}=\varphi_{T/S}(\Defres_{T/S}^Pf)
\end{equation}
for any $p$-group $P$, any $f\in F(P)$, and
any $(T,S)\in\mathcal{Y}(P)$.
We now show that the maps~$\varphi^+_P$ define a morphism
$\varphi^+:F\to \mathcal{R}_{\mathcal{Y}}G$ in~$\mathcal{F}$.\par

If $Q$ is a $p$-group and $U$ is a finite $(Q,P)$-biset, then for
each $f\in F(P)$ and each $(T,S)\in\mathcal{Y}(Q)$, one has that
\begin{eqnarray*}
\big(\mathcal{R}_{\mathcal{Y}}G(U)\circ\varphi^+_P(f)\big)_{T,S}&=&
\sum_{u\in[T\bs U/P]}(S\bs Tu)\cdot \varphi_{T^u/S^u}(\Defres_{T^u/
S^u}^Pf)\\
&=&\sum_{u\in[T\bs U/P]}\varphi_{T/S}\big((S\bs Tu)\cdot \Defres_{T^u/
S^u}^Pf\big)\\
&=&\varphi_{T/S}\left(\sum_{u\in[T\bs U/P]}(S\bs Tu)\times_{T^u/S^u}
(S^u\bs P)\cdot f\right) \,.
\end{eqnarray*}

\noindent
By Lemma~\ref{facile}, there is an isomorphism of $(T/S,P)$-bisets
$$(S\bs Tu)\times_{T^u/S^u}(S^u\bs P)\cong S\bs (Tu\times_{T^u}P) \,.$$
Moreover $Tu\times_{T^u}P\cong TuP$ as $(T,P)$-bisets.
Thus
\begin{eqnarray*}
\big(\mathcal{R}_{\mathcal{Y}}G(U)\circ\varphi^+_P(f)\big)_{T,S}
&=&\varphi_{T/S}\left(\sum_{u\in[T\bs U/P]}S\bs TuP\cdot f\right)\\
&=&\varphi_{T/S}\left(\bigl(\mathop{\bigsqcup}_{u\in[T\bs U/P]}
\limits S\bs TuP\bigr)\cdot f\right)\\
&=&\varphi_{T/S}\big((S\bs U)\cdot f\big)\\
&=&\varphi_{T/S}\big((S\bs Q \times_Q U)\cdot f\big)\\
&=&\varphi_{T/S}\big(\Defres_{T/S}^Q(U\cdot f)\big)\\
&=&\big(\varphi_Q^+(U\cdot f)\big)_{T,S}\\
&=&\big(\varphi^+_QF(U)(f)\big)_{T,S} \,.
\end{eqnarray*}
This shows that $\mathcal{R}_{\mathcal{Y}}G(U)\circ\varphi^+_P=
\varphi^+_Q\circ F(U)$, so the maps $\varphi^+_P$ define
a morphism $\varphi^+:F\to \mathcal{R}_{\mathcal{Y}}G$
in~$\mathcal{F}$.\par

Conversely, suppose that $\psi: F\to \mathcal{R}_{\mathcal{Y}}G$ is a
morphism in~$\mathcal{F}$. Then $\mathcal{O}_{\mathcal{Y}}(\psi)$ is
a morphism from $\mathcal{O}_{\mathcal{Y}}F$ to $\mathcal{O}_{\mathcal
{Y}} \mathcal{R}_{\mathcal{Y}}G$, and the latter is isomorphic
to~$G$. It is easy to check that the corresponding morphism
$\psi^-:\;\mathcal{O}_{\mathcal{Y}}F\to G$ maps the element
$f\in F(Q)$, for $Q\in\mathcal{Y}$, to $\psi_Q(f)_{Q,1}$,
where we have identified $Q/1$ with~$Q$.
In other words
$$\psi^-_Q(f)=\psi_Q(f)_{Q,1}$$
for any $Q\in \mathcal{Y}$ and any $f\in F(Q)$ .\par

We have now defined morphisms
\begin{eqnarray*}
\Hom_{\mathcal{F}_{\mathcal{Y}}}\big(\mathcal{O}_{\mathcal{Y}}F,G\big)
\longrightarrow
\Hom_\mathcal{F}\big(F,\mathcal{R}_{\mathcal{Y}}G\big)
\,, &\qquad&
\varphi \mapsto \varphi^+ \\
\Hom_\mathcal{F}\big(F,\mathcal{R}_{\mathcal{Y}}G\big)
\longrightarrow
\Hom_{\mathcal{F}_{\mathcal{Y}}}\big(\mathcal{O}_{\mathcal{Y}}F,G\big)
\,, &\qquad&
\psi \mapsto \psi^-
\end{eqnarray*}

\noindent
and we have to show that they are mutual inverse.
If $\varphi:\;\mathcal{O}_{\mathcal{Y}}F\to G$ is a morphism, then
for any $Q\in\mathcal{Y}$ and $f\in F(Q)$, we have
$$(\varphi^+)^-_Q(f)=\varphi^+_Q(f)_{Q,1}=\varphi_{Q/1}(\Defres_{Q/1}
^Qf)
=\varphi_Q(f)\mpoint$$
Hence $(\varphi^+)^-=\varphi$.\par

Conversely, if $\psi: F\to \mathcal{R}_{\mathcal{Y}}G$ is a morphism,
$P$ is a $p$-group, $f\in F(P)$, and $(T,S)\in \mathcal{Y}(P)$, then
\begin{eqnarray*}
(\psi^-)^+_P(f)_{T,S}&=&\psi^-_{T/S}(\Defres_{T/S}^Pf)\\
&=&(\psi_{T/S}\Defres_{T/S}^Pf)_{T/S,1}\\
&=&(\Defres_{T/S}^P\psi_P(f))_{T/S,S/S}\\
&=&(\psi_P(f))_{T/S,S/S}=\psi_P(f)_{T,S}\mvirg
\end{eqnarray*}
using Example~\ref{ex1} and the identification between $(T/S)/(S/S)$
and $T/S$.
Therefore $(\psi^-)^+=\psi$. \par

Thus we have proved that the correspondences
$\varphi\mapsto\varphi^+$ and $\psi\mapsto \psi^-$
provide mutual inverse isomorphisms
$$\Hom_{\mathcal{F}_{\mathcal{Y}}}\big(\mathcal{O}_{\mathcal{Y}}F,G\big)
\begin{array}{c}\raisebox{-1ex}{$\longleftarrow$}\\\raisebox{1ex}{$
\longrightarrow$}\end{array}
\Hom_\mathcal{F}\big(F,\mathcal{R}_{\mathcal{Y}}G\big)\mpoint$$
We leave to the reader the easy verification that these isomorphisms
are natural in~$F$ and~$G$.
\endpf

\result{Corollary of proof} \label{unit-defres}
Let $F\in \mathcal{F}$ and let $P$ be a $p$-group.
Let $\eta^{\mathcal Y}_F : F \longrightarrow
{\mathcal R}_{\mathcal Y}{\mathcal O}_{\mathcal Y}F$ be the unit
morphism associated with the adjunction of Theorem~\ref{right-adj}.
Then the evaluation $\eta^{\mathcal Y}_F(P)$ coincides with the map
$$ \prod_{(T,S)\in{\mathcal Y}(P)} \!\! \Defres^P_{T/S} :
\; F(P) \longrightarrow
{\mathcal R}_{\mathcal Y}{\mathcal O}_{\mathcal Y}F(P) =  \limproj
{(T,S)\in{\mathcal Y}(P)} F(T/S)\,.$$
\fresult

\pf
Via the adjunction, the unit morphism $\eta^{\mathcal Y}_F$
corresponds to the identity
$\Id : {\mathcal O}_{\mathcal Y}F \to {\mathcal O}_{\mathcal Y}F$,
that is,
$\eta^{\mathcal Y}_F=\Id^+$.
But for any $f\in F(P)$, $\Id_P^+(f)$ is defined by the expression~
(\ref{adjoint-defres}) in the proof above, that is,
the family $\Defres_{T/S}^Pf$, where $(T,S)$ runs through~$\mathcal{Y}
(P)$.
\endpf

\bigskip
\noindent
Serge Bouc, CNRS-LAMFA, Universit\'e de Picardie - Jules Verne,
33, rue St Leu,\\F-80039 Amiens Cedex~1, France.
\par\noindent {\tt serge.bouc@u-picardie.fr}

\medskip
\noindent
Jacques Th\'evenaz, Institut de G\'eom\'etrie, Alg\`ebre et Topologie,
EPFL, B\^atiment BCH, CH-1015 Lausanne, Switzerland.
\par\noindent {\tt Jacques.Thevenaz@epfl.ch}

\end{document}